\documentclass{article}
\usepackage{amsmath,latexsym,amssymb,amsthm,enumerate,amscd}
\theoremstyle{plain}
\newtheorem{theorem}{Theorem}[section]
\newtheorem{proposition}[theorem]{Proposition}
\newtheorem{corollary}[theorem]{Corollary}
\newtheorem{lemma}[theorem]{Lemma}
\theoremstyle{definition}
\newtheorem{definition}[theorem]{Definition}

\newtheorem{remark}[theorem]{Remark}
\newtheorem{example}[theorem]{Example}

\theoremstyle{plain}
\newtheorem{cor}[theorem]{Corollary}

\theoremstyle{definition}

\newtheorem*{uthm}{Theorem}

\numberwithin{equation}{section}
\numberwithin{table}{section} 
\newcommand{\Gor}{\rm{Gor}}

\DeclareMathOperator{\Grass}{\rm{Grass}}

\def\cha{\mathrm{char}\ }

\def\Hom{\mathrm{Hom}}

\def\Ann{\mathrm{Ann}\ }
\def\Proj{\mathrm{Proj}\ }

\def\PGOR{\mathbb{P}\mathrm{Gor}}

\def\Hilb{\mathrm{Hilb}}

\def\Soc{\mathrm{Soc}}

\def\Grass{\mathrm{Grass}}

\def\<{\left<}
\def\>{\right>}

\def\PGl{\mathrm{Pgl}}

\def\Sym{\mathrm{Sym}}

\def\ns{\footnotesize \it}

\def\Z{\mathfrak{Z}}
\def\W{\mathfrak{W}}

\def\PGL{\mathrm{Pgl}}
\def\max{\mathrm{max}}
\title{Artinian Gorenstein algebras of embedding dimension four:
  Components of $\PGOR(H)$ for $ H=(1,4,7,\ldots ,1)$}
\author{Anthony Iarrobino\\[.05in]
{\ns Department of Mathematics, Northeastern University, Boston, MA 02115, USA.
}\\[.2in]
Hema Srinivasan\\[.05in]
{\ns Department of Mathematics, University of Missouri, Columbia, Missouri,
USA.
}\\[.2in]}

\date{March 23, 2004}

\begin{document}

\maketitle
\centerline {{\it Dedicated to Wolmer Vasconcelos on the occasion of his
65th birthday}}

\begin{abstract} A Gorenstein sequence $H$ is a sequence of nonnegative
integers
$H=(1,h_1,\ldots ,h_j=1)$ symmetric about $j/2$ that occurs as the Hilbert
function in
degrees less or equal
$j$ of a standard graded Artinian Gorenstein algebra
$A=R/I$, where $R $ is a polynomial ring in $r$ variables and $
I$ is a graded ideal. The scheme $\PGOR (H)$ parametrizes all such
Gorenstein algebra quotients of $R$
having Hilbert function $H$ and it is known to be smooth when the embedding
dimension
satisfies $h_1 \le 3$.
The authors give a structure theorem for such
Gorenstein algebras of
Hilbert function
$H=(1,4,7,\ldots )$ when $R=K[w,x,y,z]$ and
$I_2\cong \langle wx,wy,wz\rangle $ ({Theorem~\ref{V,W}, \ref{ResI,J}}).
They also show that any Gorenstein sequence
$H=(1,4,a,\ldots ), a\le 7$ satisfies the condition
$\Delta H_{\le j/2}$ is an O-sequence ({Theorem~\ref{nonempty2},
\ref{aless7}}).
 Using these results, they show that
if $H=(1,4,7,h,b,\ldots ,1)$ is a Gorenstein sequence
satisfying $3h-b-17\ge 0$, then the Zariski closure
$\overline{\mathfrak{C}(H)}$ of the
subscheme
$\mathfrak{C}(H)\subset
 \PGOR(H)$ parametrizing Artinian Gorenstein quotients
$A=R/I$ with $I_2\cong \langle wx,wy,wz\rangle$ is a generically smooth
component of
$\PGOR(H)$ ({Theorem \ref{sevcompt}}).

They show that if in addition $8\le h\le 10$, then such $\PGOR(H)$ have
several irreducible components (Theorem \ref{sevcomp2}).  M. Boij
and others had given previous examples of certain $\PGOR(H)$  having
several components in embedding
dimension four or more
\cite{Bj2},\cite[Example C.38]{IK}.\par
 The proofs use properties
of minimal resolutions, the smoothness of
$\PGOR(H')$ for embedding dimension three \cite{Klj2}, and the Gotzmann
Hilbert scheme theorems
\cite{Gotz,IKl}.
\end{abstract}
\section{Introduction}\label{j1}\par
Let $R$ be the polynomial ring $R=K[x_1,\ldots x_r]$ over an algebraically
closed field $K$, and denote by
$M=(x_1,x_2,\ldots ,x_r)$ its maximal ideal.  When
$r=4$,  we let
$R=K[w,x,y,z]$  and regard it as
 the coordinate ring of the projective space $\mathbb
P^3$. Let $A=R/I$ be a standard graded
Artinian Gorenstein (GA) algebra, quotient of $R$. We will denote
by
$\Soc(A)=(0:M$) the socle of $A$,
the one-dimensional subvector space of $A$ annihilated by multiplication by
$M$. It is the minimal
non-zero ideal of $A$. Its degree is the \emph{socle degree} $j(A):  \,
j(A)=\max\{i\mid
A_i\not=0\}$. A sequence
$H=(h_0,\ldots ,h_j)=(1,r,\ldots , r,1)$ of positive integers symmetric about
$j/2$ is called a {\it Gorenstein sequence} of socle degree $j$, if it
occurs as the Hilbert
function of some graded Artinian Gorenstein (GA) algebra $A = R/I$. We let
$\Delta
H_i=h_i-h_{i-1}$, and denote by $H_{\le d}$ the subsequence $(1,h_1,\ldots
,h_d)$. The graded Betti
numbers of an algebra are the dimensions of the various graded pieces that
occur in the minimal graded R-resolution of
$A$.

When $r=2$, F.~H.~S. Macaulay had shown \cite{Mac1} that an Artinian
Gorenstein quotient of
$R$ is a complete intersection quotient $A=R/(f,g)$; thus, for $A$ graded,
the Gorenstein sequence must have the
form $H(A)=H(s)=(1,2,\ldots ,s-1,s,s,\ldots ,2,1)$. Also, when $r=2$ the
family $\PGOR(H(s))$ parametrizing such
Artinian quotients is smooth; its closure
$\overline{\PGOR(H(s))}=\bigcup_{t\le s}\PGOR(H(t))$ is
naturally isomorphic to the secant variety of a rational normal curve, so
is well understood (see, for example,
\cite[\S 1.3]{IK}). \par For Artinian Gorenstein algebras
$A$ of embedding dimension three ($r=3$), the Gorenstein sequences
$H(A)$, and the possible sequences $\beta$ of graded Betti numbers for $A$
given the Hilbert function $H(A)$ had
been known for some time
\cite{BE,St,D,HeTV}, see also \cite[Chapter 4]{IK}. More recently, the
irreducibility and smoothness of the
family
$\PGOR(H)$ parametrizing such GA quotients having Hilbert function $H$ was
shown by S.~J.~Diesel and
J.-O.~Kleppe, respectively
(\cite{D,Klj2}). When $r=3$, there are also several dimension formulas for
the family $\PGOR(H)$,
due to A.~Conca and G.~Valla, J.-O.~Kleppe, Y.~Cho and B.~Jung
\cite{CV,Klj2,CJ}
(see also \cite[Section 4.4]{IK} for a survey); also, M. Boij has found the
dimension of the
subfamily $\PGOR(H,\beta)$ parametrizing $A$ with a given sequence $\beta$
of graded Betti numbers
\cite{Bj3}. The closure $\overline{\PGOR(H)}$ is in general less well
understood when $r=3$, but see \cite[Theorem
5.71,\S 7.1-7.2]{IK}.\par
For embedding dimensions five or greater, it is known that a Gorenstein
sequence may be
non-unimodal: that is, it may have several maxima separated by a smaller
local minimum (\cite{BeI,BjL}).
\par
When the embedding dimension is four, it is not known whether
Gorenstein
sequences must satisfy the
condition that the first difference
$\Delta H_{\le j/2}$ is an $O$-sequence --- a sequence admissible for the
Hilbert function of some ideal of
embedding dimension three (see Definition \ref{Macexp}). Nor do we
know whether height four
Gorenstein sequences are unimodal, a weaker restriction. Little was
known about the parameter
scheme
$\PGOR(H)$ when $r=4$, except that for suitable Gorenstein sequences $H$,
it may have several irreducible
components
\cite{Bj2},\cite[Example C.38]{IK}.  We had the following questions, that
guided this portion of
our study.
\smallskip\par
$\bullet$ Can we find insight into the open problem of whether height four
Gorenstein sequences $H$ must
satisfy the condition, $\Delta H_{\le j/2}$ is an $O$-sequence?
\smallskip\par
$\bullet$ Do most schemes $\PGOR(H)$ when $r=4$ have
several irreducible components,  or is this a rare phenomenon?

\par\smallskip We now outline our main results.
We consider Hilbert sequences $H = (1, 4, 7, \cdots 1)$.  Thus, $I$ is always
a graded height four Gorenstein ideal in $K[w,x,y,z]$ whose minimal sets of
generators
include exactly three quadrics.
 First, in Theorem
\ref{V,W}, we obtain a structure theorem for Artinian Gorenstein
quotients $A=R/I$  with Hilbert function $H(A)=H$ and
with $I_2\cong \langle
wx,wy,wz\rangle$.  The proof relies on the connection between $I$ and
the intersection $J=I\cap K[x,y,z]$, which is a height three Gorenstein ideal.
  We also construct the minimal resolution of $A$ in Theorem \ref{ResI,J}.
  This allows us to determine the tangent space ${\Hom}_0(I,R/I)$ to
$A$ on $\PGOR (H)$, and to
show that under a simple condition on $H$, if such an algebra $A$ is
general enough, then $A$ is parametrized by a
smooth point of
$\PGOR(H)$ (Theorem \ref{sevcomp}).\par
 We then study the intriguing case $A=R/I$  where
$I_2\cong \langle
w^2,wx,wy\rangle$ and exhibit  a subtle connection between $A$ and a height
three
Gorenstein algebra.  We determine in Theorem~\ref{hfw2} that
the possible Hilbert functions $H=H(A)$ for such Artinian algebras $A$ satisfy
\begin{equation}\label{simpleH}
H=H'+(0,1,\ldots ,1,0)
\end{equation}
where $H'$ is a height three Gorenstein sequence. \par
 Our result pertaining to the first question is
\begin{uthm} (Theorem \ref{nonempty2}, Corollary
\ref{nonempty3}, Proposition \ref{aless7})
All Gorenstein sequences of the form $H=(1,4,a,\ldots ),\, a\le 7$  must
satisfy the condition
that $\Delta H_{\le j/2}$ is an
$O$-sequence.
\end{uthm}
To show this we eliminate
potential sequences not satisfying
the condition by frequently using the symmetry of the minimal resolution of
a graded Artinian Gorenstein algebra $A$, the Macaulay bounds on the
Hilbert function,
and the Gotzmann Persistence and Hilbert scheme theorems (Theorem
\ref{MacGo}). However, these methods
 do not
extend to all height four Gorenstein sequences, and we conjecture that not
all will satisfy the
condition that $\Delta H_{\le j/2}$ is an
$O$-sequence (see Remark~\ref{SIconj}).
\par
We then combine these results with a well known construction of
Gorenstein ideals from sets of
points to obtain our theorem concerning irreducible components of
$\PGOR (H)$

\begin{uthm}(Theorem
\ref{sevcomp2}\ref{sevcomp2A}) Let $H=(1,4,7,h,b,\ldots ,1)$ be a
Gorenstein sequence
satisfying $8\le h\le 10$ and
$3h-b-17\ge 0$. Then $\PGOR (H)$ has at least two components. The first
is the Zariski
closure of the subscheme $\mathfrak C(H)$
of $
 \PGOR(H)$ parametrizing Artinian Gorenstein quotients
$A=R/I$ for which $I_2$ is $\PGL(3)$-isomorphic to $ \langle
wx,wy,wz\rangle$. The second
component parametrizes quotients of the coordinate rings of certain
punctual schemes in $\mathbb
P^3$.
\end{uthm}

\section {Notation and basic results}

In this Section we give definitions and some basic
results that we will need. Recall that $R=K[w,x,y,z]$ is the
polynomial ring with the standard grading over an algebraically closed
field, and that we consider only graded ideals $I$.
\par
 Let $V\subset R_v$ be a vector subspace. For
$u\le v$ we let
$V:R_u=\langle f\in R_{v-u}\mid R_u\cdot f\subset V\rangle$.
We state as a lemma a result of Macaulay \cite[Section
60ff]{Mac1} that we will use frequently.
\begin{lemma}\label{MacD}(F.H.S. Macaulay). Let $\cha K=0$ or $\cha K>j$.
There is a one-to-one
correspondence between graded
Artinian Gorenstein algebra quotients
$A=R/I$ of
$R$ having socle degree
$j$, on the one hand, and on the other hand, elements
$F\in \mathcal R_j$ modulo $ K^{\ast}-$action. where $\mathcal
R=K[W,X,Y,Z]$, the dual polynomial
ring.
 The correspondence is given by
\begin{align}\label{eMac1}
 I&=\Ann F = \{h\in R\mid h\circ F=h(\partial /\partial W,\ldots ,\partial
/\partial
Z)\circ F=0\};\notag\\
F&=(I_j)^{\perp}\in \mathcal R_j \mod K^{\ast} .
\end{align}
Here $F$ is also the generator of the $R$-submodule $I^\perp\subset
\mathcal R, I^\perp=\{ G\in
\mathcal R\mid h\circ G=0 \text { for all } h\in I\}$. The Hilbert function
$H(R/I)$ satisfies
\begin{equation}
H(R/I)_i=\dim_K (R\circ F)_i=H(R/I)_{j-i}.
\end{equation}
Furthermore, for $i\le j$, $I_i$ is determined by $I_j$ or by $F$ as follows:
\begin{equation}\label{ancclosure}
I_i=I_j:R_{j-i}=\{h\in R_i\mid h\cdot R_{j-i}\subset I_j\}=\{ h\in R_i\mid
h\circ(R_{j-i}\circ F)=0\}.
\end{equation}
When $\cha K=p>j$ the statements are analogous, but we must replace
$K[W,X,Y,Z]$ by the ring of
divided powers
$\mathcal D$, and the action of $R$ on $\mathcal D$ by the contraction action
(see below).
\end{lemma}
\begin{proof} For a modern proof see \cite[Lemmas 2.15, 2.17]{IK}. For a
discussion of the use of
the divided power ring when $\cha K = p$ see also \cite[Appendix A]{IK}.
\end{proof}
\begin{corollary}\label{include} Let $A=R/I$ be a graded Artinian
Gorenstein algebra
of socle degree $j$. Let $J=I_\Z$ be a saturated ideal defining a
scheme $\Z\subset \mathbb P^{3}$,
such that for some $i, 2\le i\le j$, $\Z=\Proj(R/(J_i))$, with
$J_i\subset I_i$. Then for $0\le u\le i$ we have $J_u\subset I_u$. If also
$J_i=I_i$, then for such $u$,
$J_u=I_u$.
\end{corollary}
\begin{proof} Let $0\le u\le i$. Since
$J$ is its own saturation, we have $J_u=J_k:R_{k-u}$ for large $k$, so we have
\begin{equation*}
J_u=J_k:R_{k-u}=\{ J_k:R_{k-i}\} : R_{i-u}=J_i:R_{i-u}.
\end{equation*}
Now \eqref{ancclosure} implies that for $0\le u\le i$
\begin{equation*}
I_u=I_j:R_{j-u}=\{ I_j:R_{j-i}\} :R_{i-u}=I_i:R_{i-u}.
\end{equation*}
This completes the proof of the relation between $I_\Z$ and $I$.
\end{proof}\par
Note that \cite[Example 3.8]{I}, due to D. Berman, shows that one cannot
conclude that $J\subset I$
in
Corollary~\ref{include}. For let $I=(x^3,y^3,z^3)$, and
let $J$ be the saturated ideal
$J=(x^2y^3,y^2z^3,x^3z^2,x^2y^2z^2)$, a local complete intersection of
degree 18 defining a
punctual scheme concentrated at the points $(1,0,0),(0,1,0)$ and $(0,0,1)$.
Then we have
$J_5\subset I_5$ but $x^2y^2z^2\in J$ so $J\nsubseteq I$.
\par
 We suppose $R=K[w,x,y,z]$.  Let ${ \mathcal D}=K_{DP}[W,X,Y,Z]$ denote the
divided power
algebra associated  to $R$: the basis of $\mathcal D_j$ is $\{ W^{[j_1]}\cdot
X^{[j_2]}\cdot Y^{[j_3]}\cdot Z^{[j_4]}, \sum j_i=j\}$. We let $x^i\circ
X^{[j]}
= X^{[j-i]}$ when
$j\ge i$ and zero otherwise;  this action extends in a natural way to the
contraction action of
$R$ on
$\mathcal D$. Multiplication in $\mathcal D$ is determined by $X^{[u]}\cdot
X^{[v]}=\binom{u+v}{v}X^{[u+v]}$. By $(\alpha X+ Y)^{[u]}, \alpha\in K$ we
mean $\sum_{0\le
i\le u}\alpha^iX^{[i]}\cdot Y^{[u-i]}$: this is $(\alpha X+Y)^u/u!$ when
the latter makes sense.
When
$\cha K=0$, or $\cha K>j$ we may replace
$\mathcal D$ by the polynomial ring
$\mathcal R=K[W,X,Y,Z]$ with $R$ acting on $\mathcal R$ as partial
differential operators
\eqref{eMac1}, and we replace all $X^{[u]}$ by $X^u$, and $(\alpha
X+Y)^{[u]}$ by $(\alpha
X+Y)^u$. \par The inverse system
$I^\perp\subset
\mathcal D$ of the ideal
$I\subset R$ satisfies
\begin{equation}
I^{\perp} = \{G\in  K_{DP}[W,X,Y,Z], h\circ G = 0 \text { for all } h\in I\},
\end{equation}
 and it is an $R$-submodule of $\mathcal D$ isomorphic to the dual module
of $A=R/I$.
 When $A=R/I$ is graded Gorenstein of socle degree $j$, then by Macaulay's
Lemma \ref{MacD} the
inverse system is principal, generated by $F\in \mathcal D_j$: we call $F$
the \emph{dual
generator} of $A$ or for $I$. Thus, we may parametrize the algebra $A$ by
the class of
$F\mod$ nonzero $ K^\ast -$multiple, an element of the projective space
$\mathbb
P^{N-1},N={\binom{j+3}{j}}$. Given a Gorenstein sequence
$H$ of socle degree $j$ (so $H_j\not= 0, H_{j+1}=0$) we let $\PGOR(H)\subset
\mathbb P^{N-1}$ denote  the scheme parametrizing the
family of all GA quotients
$A=R/I$ having Hilbert function $H$. Here, we use the scheme structure given by
the catalecticants, and described in
\cite[Definition 1.10]{IK}.   A
``geometric point'' $p_A$ of $\PGOR(H)$
parametrizes a Artinian Gorenstein quotient $A=R/I$ of $R$ having Hilbert
function $H$.

We now state Macaulay's theorem characterizing Hilbert functions or
$O$-sequences, and the
version of the Persistence
and Hilbert Scheme theorems of G. Gotzmann that we will use \cite{Gotz}.

  Let $d$ be a positive integer. The $d$-th
Macaulay coefficients of a positive
integer
$c$ are the unique decreasing sequence of
non-negative integers $k(d),\ldots ,k(1)$ satisfying
\begin{equation*}
c=\binom{k(d)}{d}+\binom{k(d-1)}{d-1}+\cdots
+\binom{k(1)}{1};
\end{equation*}
 We denote by $c^{(d)}$  the integer
\begin{equation}\label{eMacExp}
c^{(d)}=\binom{k(d)+1}{d+1}+\binom{k(d-1)+1}{d}+\cdots +\binom{k(1)+1}{2}.
\end{equation}

Then,
 the Hilbert polynomial
$p_{c,d}(t)$ for quotients $B$ of the polynomial ring $R$, such that $B$ is
regular in degree
$d$ and $H(B)_d=c$
satisfies
\begin{equation}\label{eMacExp2}
p_{c,d}(t)=\binom{k(d)+t-d}{t}+\binom{k(d-1)+t-d}{t-1}+\cdots
+\binom{k(1)+t-d}{t-d}.
\end{equation}
The length of the
$d$-th Macaulay expansion of $c$, or of the Macaulay expansion of the
polynomial $p_{c,d}$, is the number of
$\{ k(i)\mid k(i)\ge i\} $, equivalently, the number of nonzero binomial
coefficients in the Macaulay
expansion, and this is well known to be the Gotzmann regularity degree of
$p_{c,d}$ (\cite[Theorem
4.3.2]{BH}).
\begin{theorem}\label{MacGo} Suppose that $1\le c\le \dim_kR_d$, and $I$ is
a graded ideal of
$R=K[x_1,\ldots x_r]$.
\begin{enumerate}[i.]
\item\label{Macgrowth}\cite{Mac2} If $H(R/I)_d=c$, then $H(R/I)_{d+1}\le
c^{(d)}$ (Macaulay's inequality).
\item\label{Gotzper}\cite{Gotz} If $H(R/I)_d=c$ and $H(R/I)_{d+1}=c^{(d)}$,
then
$\Proj(R/(I_d))$ is a projective scheme in $\mathbb P^{r-1}$ of Hilbert
polynomial $p_{c,d}(t).$

In particular
 $H(R/(I_d))_{k}=p_{c,d}(k)$ for $k\ge d$, and $H'=H(R/(I_d))$ has extremal
growth ($h'_{k+1}={{h'}_k}^{(k)}$) in each degree
$k$ to
$k+1, k\ge d$.
\end{enumerate}
\end{theorem}
\begin{proof} For a proof of Theorem \ref{MacGo}\eqref{Macgrowth} see
\cite[Theorem 4.2.10]{BH}.
For a proof of the persistence (second) part of Theorem \ref{MacGo}
\eqref{Gotzper} see
\cite[Theorem 4.3.3]{BH}; for the Gotzmann Hilbert scheme theorem see
\cite[Satz1]{Gotz}, or the
discussion of
\cite[Theorem C.29]{IKl}.
\end{proof}
\begin{definition}\label{Macexp}
A sequence of nonnegative integers $H=(1,h_1,\ldots ,h_d,\ldots ),$ is said
to be an
$O$-{\emph{sequence}}, or to be \emph{admissible} if it
satisfies Macualay's inequality of Theorem \ref{MacGo}\eqref{Macgrowth}
for each integer $d\ge 1$.
\end{definition}
Recall that the regularity degree $\sigma(p)$ of a Hilbert polynomial
$p=p(t)$ is the smallest degree for
which all projective schemes $\Z$ of Hilbert polynomial $p$ are
Castelnuovo-Mumford regular in  degree
less or equal
$\sigma (p)$.  G. Gotzmann and D. Bayer showed that this bound is the
length $\sigma(p)$ of the
Macaulay expression for $p$ \cite{Gotz,Ba}: for an exposition
and proof see \cite[Theorem 4.3.2]{BH}; also see
\cite[Definition C.12, Proposition C.24]{IKl}, which includes some
historical remarks. As an easy consequence we
have
\begin{corollary}\label{regdeg} The regularity degree of the polynomial
$p(t)=at+1-\binom{a-1}{2}+b$ where $a >0,b\ge 0$ satisfies $
\sigma (p)=a+b$. These
Hilbert polynomials cannot occur with
$b<0$. In particular we have, the regularity degree of the polynomial
$p(t)=3t+b, b\ge 0$ is $3+b$,
of $p(t)=2t+1+b, b\ge 0$ is $b+2$, and of $p(t)=t+1+b,b\ge 0$ is $b+1$. The
regularity of
the constant polynomial $p(t)=b$ is $b$.
\end{corollary}
\proof One has for $p(t)=at+1-{a-1\choose 2}+b$, the following sum,
equivalent to a
Macaulay  expansion as in \eqref{eMacExp2} of length $a+b$,
\begin{equation*}
\begin{split}
p(t)={t+1\choose
1}&+{t+1-1\choose 1}+{t+1-2\choose 1}+\dots + {t+1-(a-1)\choose
1}+\\
&{t-a\choose 0}+{t-(a+1)\choose 0}+\dots + {t-(a+b-1)\choose 0}.
\end{split}
\end{equation*}
\par
\begin{corollary}\label{oseq}
Let $H$ be a Gorenstein sequence of socle degree $j$, and suppose for that
some $d< j$,
$h_{d+1}=(h_{d})^{(d)}$ is extremal in the sense of Theorem
\ref{MacGo}\eqref{Macgrowth}. Then
$\Delta H_{\le d+1}$ is an $O$-sequence.
\end{corollary}
\begin{proof} Theorem \ref{MacGo}\eqref{Gotzper} and Corollary
\ref{include} show the existence of a scheme
$\Z\subset
\mathbb P^{r-1}$ satisfying $h_u=H(R/I_\Z)_u$ for $u\le d+1$. Since
$I_\Z$ is saturated and thus $R/I_\Z$ has depth at least one, there is a
homogeneous degree one
nonzero divisor, implying that the first difference $\Delta (H(R/I_\Z))$ is
an $O$-sequence.
\end{proof}\par
\begin{remark}
The assertion of Corollary \ref{oseq} as well as those of Corollary
\ref{include} are valid more
generally for graded Artinian algebras having socle only in degree $j$
(level algebras), or those
having socle only in degrees greater or equal $j$.
\end{remark}
As an example of the application of Theorem \ref{MacGo}, we determine below
the Gorenstein
sequences
\linebreak $H=(1,4,7,h,7,4,1)$ that
 occur, having socle degree 6.
\begin{corollary}\label{nonempty1} The sequence $H=(1,4,7,h,7,4,1)$ is a
Gorenstein sequence
if and only if
$7\le h\le 11$.
\end{corollary}
\begin{proof} From Macaulay's extremality Theorem
\ref{MacGo}\eqref{Macgrowth} we have
$H(3)=h\le H(2)^{(2)}=7^{(2)}=11,$ and
$H(4)=7\le h^{(3)}$ which implies $h\ge 6.$
 Now $H=(1,4,7,6,7,4,1)$
implies that
the growth of $H_3=6$ to $H_4=7$ is maximum, since
$6=\binom{4}{3}+\binom{2}{2}+\binom{1}{1}$, while
$7=6^{(3)}=\binom{5}{4}+\binom{3}{3}+\binom{2}{2}$. Corollary~\ref{oseq}
shows this is impossible.
\end{proof}\par
For a subscheme $\Z\subset \mathbb P^3$ we will denote by
$H_\Z=H(R/I_\Z)$ its Hilbert function, sometimes
called its postulation; here $I_\Z\subset R$ is the saturated ideal
defining $\Z$. Inequalities among Hilbert
functions are termwise. The following result is well-known and easy to
show, since a degree-d punctual scheme can
cut out at most d conditions in a given degree.
\begin{lemma}\label{upbd} Let $\Z=W\cup \Z_1\subset \mathbb P^3$, be a
subscheme of $\mathbb P^3$,  where $W$
is  a degree
$d$ punctual scheme. Then for all $i$, $(H_\Z)_i \le (H_{\Z_1})_i+d$.
\end{lemma}
\begin{proof} We have (the first inequality is from Maroscia's result
\cite{Mar}, see
\cite[Theorem 5.1A]{IK})
\begin{equation}
\begin{split}
d\ge (H_W)_i=&\dim R_i-\dim (I_W)_i\ge \dim (I_{\Z_1})_i-\dim (I_W\cap
I_{\Z_1})_i\\=&H(R/(I_W\cap I_{\Z_1}))_i-H(R/I_{\Z_1})_i\ge
(H_\Z)_i-(H_{\Z_1})_i.
\end{split}
\end{equation}
\end{proof}
\section{Nets of quadrics in $\mathbb P^3$, and Gorenstein
ideals}\label{netsec}

In Section \ref{Netsquad} we give preparatory material on nets of quadrics,
and on the Hilbert schemes
of low degree curves in $\mathbb P^3$. In Section \ref{struc} we prove a
structure theorem for Artinian
Gorenstein algebras
$A=R/I$ of Hilbert function
$H(A)=(1,4,7,\ldots )$ for which the net of quadrics $I_2$ has a common factor
and is isomorphic
after a change of variables to $\langle wx,wy,wz\rangle$ (Theorem
\ref{V,W}). We then determine the dimension of the tangent space to $\PGOR
(H)$ at a point
parametrizing such an ideal; we also show that when $H$ has socle degree 6,
the subfamily parametrizing such Gorenstein algebras is an irreducible
component of $\PGOR(H)$
(Theorem \ref{sevcomp}), a result which we will later generalize to
arbitrary socle degree
(Theorem \ref{sevcompt}). In Section
\ref{w2sec} we determine  the possible Hilbert functions $H(A),A=R/I$ when
$I_2=\langle w^2,wx,wy
\rangle$ (Theorem
\ref{hfw2}).
\subsection{Nets of quadrics}\label{Netsquad}
  Three homogeneous quadratic polynomials $f,g,h$ in $R=K[w,x,y,z]$ form a
family
$\alpha_1f+\alpha_2g+\alpha_3 h,
\alpha_i\in K$, comprising a net of quadrics in $\mathbb P^3$. Here we will
use the term net also
for the vector space span $V=\langle f,g,h\rangle$. We divide these
families according to the
number of linear relations among the three quadrics.  We now show that
they can have at most 3 linear
relations.  Let
$(I_2) = (f,g,h)$ be  the ideal generated by a net of quadrics $I_2=\langle
f,g,h\rangle$.  Then
$H(R/(I_2)) = ( 1,4,7, h,
\cdots)$, where
$h\leq 11=7^{(3)} $ by Macaualay's growth condition.  When there are no
relations
$H(R/(I_2))_3=20-12=8$, so the number of linear relations on the net of
quadrics
$\langle f,g,h\rangle$ is no greater than $ 11-8 =3$, as claimed.  \par
  Nets of quadrics in $\mathbb P^3$ have been extensively
studied geometrically, earlier by
W.~L.~Edge and others, more recently  by C.T.C. Wall and others for their
connections with mapping
germs, and instantons. I. Vainsecher and also G. Ellingsrud, R. Peine, and
S.A.
Str\o{m}me have showed that the Hilbert scheme of twisted cubics in $\mathbb
P^3$ is a blow-up of the
family $\mathfrak {F}_{\mathrm{RNC}}$ of nets of quadrics arising as minors
of a
$2\times 3$ matrix (Definition \ref{defF}) along the
sublocus of those nets having a common factor. Nets of quadrics are
parametrized by the
Grassmanian $\mathbb{G}=\Grass (3,R_2)\cong \Grass(3,10)$, of dimension 21. It
is easy to see that up to
isomorphism under the natural $\PGl(3)$ action,  the vector spaces
$V=\langle f,g,h\rangle\subset R_2$ have a 6 dimensional family of orbits,
as $\dim \Grass(3,10)-\dim
\PGl(3)=21-15=6$, and the stabilizer of a general enough net is
finite.  In this section, we determine the irreducible components of the
subfamily $\mathfrak{F}$ of
nets having at least one linear relation (Lemma \ref{netsq0}), and also the
possible graded Betti
numbers for the algebras $R/(V)$, for nets $V\in \mathfrak{F}$ (Lemma
\ref{netsq}).

\begin{definition}\label{defF} We denote by $\mathfrak{F}\subset\mathbb{G}=
\Grass(3,R_2)$ the subfamily of nets of
quadrics, vector spaces
$V=\langle f,g,h\rangle\subset R_2$, for which $f,g,h$ have at least one
linear relation
\begin{equation}\label{linrel}
\alpha_1f+\alpha_2g+\alpha_3h=0,
\exists\alpha_i\in R_1=\langle w,x,y,z\rangle.
\end{equation}  \par
We denote by $\mathfrak{F}_i\subset\mathbb{G}=\Grass(3,R_2)$ the subfamily of
$\mathfrak{F}$ consisting of
those nets that have exactly $i$ linear relations, $i=1,2,3$.  We denote
by
$\mathfrak{F}_{\mathrm{RNC}}\subset
\mathfrak{F}_2$ the subset of nets defining twisted cubic curves, and by
$\mathfrak{F}_{sp}$ the
subset of nets $\PGL(3)$ isomorphic to $\langle w^2,wx,wy\rangle$.
\end{definition}
\begin{lemma}\label{F1} The family $\mathfrak{F}_1$ comprises those nets
that can be written
$V=\langle \ell\cdot U,h\rangle$, where $\ell\in R_1$ is a linear form,
$U\subset R_1$ is a two
dimensional subspace of linear forms, and $h$ is not divisible by either
$\ell$ or by any element
of $U$.\par
Up to isomorphism
$V\in
\mathfrak{F}_1$ may be
written either
$V=\langle xw,yw,h\rangle$ for some quadric
$h$ divisible neither by $w$ nor by any element of $\langle x,y\rangle$,
or $V=\langle
w^2,wx,h\rangle$ with
$h$ divisible by no element of $\langle w,x\rangle$.
\end{lemma} 
\begin{proof} First consider nets $V=\langle f,g,h\rangle $ having no
two dimensional subspace with a common factor: we show that $V$ cannot be
in $\mathfrak {F}_1$. When the coefficients of a relation as in
\eqref{linrel} form an $m$-sequence, a simple argument given in the proof
of Lemma \ref{netsq0} shows that $V\in
\mathfrak F_2$, and is determinantal (see equation \eqref{rnc}).\par
 Now assume that $V$ has a relation as in 
\eqref{linrel} such that $\dim_K\langle
\alpha_1,\alpha_2,\alpha_3\rangle =2$; after a change of basis in $R$ we
may suppose that
$xf+yg+(x+y)h=0$ where $h$ may be zero.
Replacing $f$ by $f+h$, and $g$ by $g+h$, we obtain
$xf=-yg$. Thus $V$ may be written $V=\langle U\ell,h\rangle$ with
$\ell=f/y$ and
$U=\langle x,y\rangle$, and, evidently if $V\in \mathfrak {F}_1$ then
$h$ is not divisible by $\ell$ nor by any element of $U$.  We have shown
the first claim of the Lemma. The second follows.
\end{proof}\par
 As we shall see
below,
$\mathfrak {F}_2$ has
 $\mathfrak{F}_{\mathrm{RNC}}$ as open dense subset.
Evidently, the family ${\mathfrak{F}}_3$ of nets $V$ having a common
factor, contains as open
dense subset the ${\PGL}(3)$ orbit of
$V=\langle wx,wy,wz\rangle$; the family also contains $\mathfrak
{F}_{sp}$, the orbit of $\langle
w^2,wx,wy\rangle$. \par
The dimension calculations of the following lemmas are elementary;
recall that $\dim\mathbb{G}=21$. The results about closures also involve standard
methods but are more subtle: for example to identify
$\mathfrak{F}_{sp}$ with $\overline{\mathfrak
{F}_2}\cap\mathfrak{F}_3$ we rely on previous work on the
closure of the family of rational normal curves, such as \cite{No,PS,Va,Lee}.
\begin{lemma}\label{netsq0}{\sc Components of $\mathfrak{F}$}: The
subfamily
$\mathfrak {F}\subset\mathbb{G}=\Grass(3,R_2)$ 
parametrizing quadrics having at least one linear relation, has two
irreducible components,
$\overline{\mathfrak {F}_1}$ and $\overline{\mathfrak
{F}_2}=\overline{\mathfrak{F}_{RNC}}$, of codimensions 7  and 9,
respectively in
$G$. They satisfy
\begin{enumerate}[i.]
 \item\label{netsq0i} The intersection
$\overline{\mathfrak {F}_1}\cap
\mathfrak {F}_2$, has an open dense subset parametrizing nets
isomorphic to $\langle
wx,wy,xz\rangle$; this intersection has codimension 11 in
$G$.
\item\label{netsq0ii} We have
$\overline{\mathfrak{F}_1}-\mathfrak{F}_1=(\overline{\mathfrak{F}_1}\cap
\mathfrak{F}_2)\cup
\mathfrak{F}_3$. Each element of $\mathfrak{F}_2$ has a basis consisting of
minors of a $2\times 3$ matrix of linear forms.  
\item\label{netsq0iii} The
locus
$\mathfrak{F}_3\subset
 \overline{\mathfrak{F}_1}$ has codimension 15 in $G$;
$\mathfrak{F}_3-\mathfrak{F}_{sp}$ consists
of nets isomorphic to $\langle wx, wy, wz\rangle$. The locus
$\mathfrak{F}_{sp}=\overline{\mathfrak{{F}_2}}\cap\mathfrak{F}_3$, and is a
subfamily of
codimension 16 in $G$.
\end{enumerate}
\end{lemma}
\begin{proof} 
We first calculate $\dim \mathfrak{F}_1$. By Lemma \ref{F1} $V\in
\mathfrak{F}_1$ may
be written as $\langle\ell\cdot U,h\rangle$, where $\ell\in R_1$ and $\cdot
U\subset R_1$ is a two dimensional subspace, and $h$ is not divisible
by $\ell$ nor by any element of $U$. Since there is a single
linear relation, $V$ determines both
$\ell$ and $U$ uniquely. Thus, 
there is a surjective morphism
\begin{equation*}
\pi_1: \mathfrak{F}_1\to \mathbb P^3\times \Grass(2,R_1):
\pi_1(V)=(\ell,U),
\end{equation*}
The fibre of $\pi_1$ over the pair $(\ell, U )$ corresponds to the
choice of $h$; given $V$, $h$ is
unique up to constant multiple, mod an element of $\ell\cdot U$. Thus, the
fibre of $\pi_1$ is
parametrized by an open dense subset of the projective space  $\mathbb
P(R_2/\langle \ell\cdot
U\rangle )$, of dimension 7. Thus,
$
\mathfrak{F}_1$ has dimension 14, and  codimension 7 in $G$.\par
We next show that $\mathfrak {F}_2$ contains $\mathfrak{F}_{RNC}$ as
dense open subset. When there is a linear relation for $V$ as in
\eqref{linrel} whose coefficients
$\alpha_i$ are a length 3 regular sequence we may suppose after a
coordinate change that
$xf+yg+zh=0$; letting $f=uz+f_1, g=vz+g_1$, with $f_1,g_1$ relatively prime
to $z$, we obtain $h=-(ux+vy)$,
and
$xf_1=-yg_1$, whence there is a linear form $\beta\in R_1$ with
$f=uz+y\beta, g=vz-x\beta$, and $(f,g,h)$ is
the ideal of $2\times 2$ minors of
\begin{equation}\label{rnc}
\begin{pmatrix}
u&v&\beta\\
-y&x&z
\end{pmatrix}.
\end{equation}
When also $(f,g,h)$ has height two,
then
$V$ is an element of
$\mathfrak{F}_2$ determining a twisted cubic in $\mathbb P^3$; for a dense
open subset of
such elements of $\mathfrak{F}_2$ one may up to isomorphism choose in
\eqref{rnc} the triple
$(u,v,\beta)=(x,z,w)$. Otherwise,if $f,g,h$ is not Cohen-Macualay of height two,
$V$ has a common linear factor, and it is well known that then $V\in
\mathfrak{F}_{sp}=
\overline{\mathfrak{F}_2}\cap\mathfrak{F}_3$ \cite{Lee,PS,Va}.\par
We now consider those nets $V\subset \mathfrak {F}_2$ for which there is
no linear relation as in \eqref{linrel} whose coefficients form a 
length three $m$-sequence. By the proof of Lemma \ref{F1} such a net has
the form $V=\langle Uáw, h\rangle$, with $U\subset R_1$, and it
thus lies in the closure of $\mathfrak {F}_1$. It is easy to see that the
most general element of
$\overline{\mathfrak{F}_1}\cap
\mathfrak{F}_2$ is a net isomorphic to $\langle wx,wy,xz\rangle$:  for
when
$V=\langle wx,wy,h\rangle$
has a second linear relation, either $w$ divides $h$ and $V\in
\mathfrak{F_3}$, or some $ax+by$
divides
$h$, and after a change in basis for $R_1$,  $V\cong \langle
wx,wy,xz\rangle$. A similar
discussion for $\langle w^2,wx,h\rangle$ completes the proof that any
element of
$\overline{\mathfrak{F}_1}\cap
\mathfrak{F}_2$ is in the closure of the orbit of $V=\langle
wx,wy,xz\rangle$, which is also the
determinantal ideal of $\begin{pmatrix}
x+y&y&0\\
z&z&w
\end{pmatrix}$.  This shows 
also that
$\overline{\mathfrak{F}_1}\cap
\mathfrak{F}_2\subset \overline{\mathfrak{F}_{\mathrm{RNC}}}$, and
completes the proof that $\mathfrak {F}_2$ contains
$\mathfrak{F}_{\mathrm{RNC}}$ as dense open subset.\par We recall that
$\dim
\mathfrak{F}_{\mathrm{RNC}} =12$. A twisted cubic -- a rational normal
curve of degree three -- is determined by the choice of four
degree three forms in the polynomial ring $K[x,y]$, up to common
$K^\ast$-multiple, mod the action
of $\mathrm{{\PGL}} (1)$, yielding dimension $4\cdot 4-4=12$
\cite{PS}.\par We have that
$\overline{\mathfrak{F}_1}$ and
$\overline{\mathfrak{F}_2}$ define two distinct
irreducible components of $\mathfrak{F}$, since the subfamily
$\mathfrak{F}_2$ parametrizing nets for which there are two linear
relations, cannot specialize
to any net
$V=\langle f,g,h\rangle$ for which $f,g,h$  have a single linear relation;
and $\mathfrak{F}_1$, parametrizing
nets $V$ each containing a subspace of the form $\ell \cdot U$, cannot specialize
to a vector space $V$ for which
the ideal
$(V)$ is the prime ideal of a twisted cubic. This completes the proof of
the initial claims of the lemma.
\par We now complete the proof of \eqref{netsq0i}, by determining the
dimension of
$\overline{\mathfrak{F}_1}\cap
\mathfrak{F}_2$, which is by the above argument equal to the dimension of
the
$\mathrm{{\PGL}(3)}$- orbit $\mathcal B$ of
$\langle wx,wy,xz\rangle$. For $W=\langle w'x',w'y',x'z'\rangle\in \mathcal
B$, the unordered pair
of linear forms
$(w',x')$, each mod $K^\ast$-multiple is uniquely determined by $W$ (as
each divides a two
dimensional subspace of $W$): thus  there is a morphism $\pi :\mathcal B\to
\Sym^2(\mathbb P^3)$, from
$\mathcal B$ to the symmetric product, whose image is the non-diagonal
pairs. Spaces $W$ in the
fibre of $\pi$ over
$(w',x')$ are determined by the choice of the two 2-dimensional subspaces,
the first $\langle
x',y'\rangle$ containing $x'$, the second $\langle w',z'\rangle$ containing
$w'$. Thus, a space $W$
in the fibre is determined by the choice of
$y'\in R_1/\langle x'\rangle$ and
$z'\in R_1/\langle w'\rangle$, each up to
$K^\ast$-multiple, and these choices are each made in an open dense subset
of $\mathbb P^2$ (as
$z'$ must not equal $x'\mod w'$ for $W\in \mathcal B$).  Thus, the fibre
$\pi^{-1}(w',x')\subset\mathcal B$ is isomorphic to an open dense subset of
$\mathbb P^2\times
\mathbb P^2$. It follows that $\mathcal B$ and $\overline{\mathfrak{F}_1}\cap
\mathfrak{F}_2$ have dimension 10, and codimension
11 in
$G$.\par We now show the claim in \eqref{netsq0ii} that
$\overline{\mathfrak{F}_1}-\mathfrak{F}_1=(\overline{\mathfrak{F}_1}\cap
\mathfrak{F}_2)\cup
\mathfrak{F}_3$. Suppose that $V\in
\overline{\mathfrak{F}_1}-\mathfrak{F}_1$; then evidently there is a
two-dimensional subspace $V_1\subset V$ having a common factor $ V_1=\elláU$.
Letting
$V=\langle V_1,h\rangle$ then $V\in \mathfrak{F}_2$ implies
$h$ must have a common divisor with an element of $V_1$. Thus, up to
$\mathrm{{\PGL}(3)}$ isomorphism we have
 $V=\langle wx,wy,xz\rangle$ or $V=\langle w^2,wx,xz\rangle$, both
in $\mathfrak{F}_2$ (we may ignore $w$ is a common factor of $V$ since then $V\in
\mathfrak {F}_3$). Each of these spaces has basis the minors of a $2\times 3$
matrix of linear forms. This with \eqref{rnc} above completes the proof of
\eqref{netsq0ii}.\par The family
$\mathfrak{F}_3$ has as open dense subset the orbit $\mathcal B'$ of
$V=\langle wx,wy,wz\rangle $.
An element $W'=w'V', V'\subset R_1$ of $\mathcal B'$ is determined by a
choice of $w'\in R_1$
and a codimension one vector space $V'\subset R_1$, thus $\mathcal B'$ is
an open in $\mathbb
P^3\times \mathbb P^3$, so has dimension six, codimension 15 in $G$. \par
The claim in \eqref{netsq0iii} that the locus
$\mathfrak{F_{sp}}=\overline{\mathfrak{F}_2}\cap
\mathfrak{F_3} $ follows from the well known
classification of the
specializations of rational normal curves \cite{PS,Lee}; the dimension
count for this locus is five, 3 for the
choice of $w$, and 2 for the choice of $\langle x,y\rangle\subset
R_2/\langle w^2\rangle$. This
completes the proof of Lemma \ref{netsq0}.
\end{proof}\par

\begin{lemma}\label{netsq}{\sc Minimal resolutions for nets of quadrics in
$\mathfrak{F}$}. There are exactly
 three possible sets of graded Betti numbers for the ideal generated by a
net of quadrics in
$\mathfrak{F}$ (those having at least one linear relation):\par
\begin{enumerate}[i.]
\item\label{netsqii} Those
$V$ in the family $\mathfrak{F}_1$ have graded Betti numbers that of
$(wx,wy,z^2)$, with a single linear and
two quadratic relations, and Hilbert function
$H=H(R/(V))=(1,4,7,9,11,13,\ldots)$ where $H_i=2i+3$ for $i\ge 2$. Such $V$
define a curve of
degree 2, genus -2. (See Lemma \ref{deg2}).
\item\label{netsqi}
 For $V\in\mathfrak{F}_2$, the ideal $(V)$ is Cohen-Macaulay
of height two, the
Hilbert  function
$H=H(R/(V))=(1,4,7,10,13,\ldots )$ where $H_i=3i+1$ for $i\ge 0$, and $V$
has the standard
determinantal minimal resolution with two linear
relations.
\item\label{netsqiii} Those $V$ in the family $\mathfrak{F}_3$ have
graded Betti numbers that of $(wx,wy,wz)$.
\end{enumerate}
\end{lemma}

\begin{proof}
For \eqref{netsqii}, Lemma \ref{F1} implies that the ideal determined by
an element $V=(wx,wy,h)$ of $\mathfrak{F}_1$ is cut out from $R/(wx,wy)$ or
$R/(w^2,wx)$ by the nonzero-divisor $h$, hence the minimal resolution of $R/(V)$
is that of
$R/(wx,wy,z ^2)$.
 For \eqref{netsqi} let $V\in \mathfrak{F}_2$.  Then by Lemma
\ref{netsq0}\eqref{netsq0ii}, $V$ is has a basis consisting of the minors of a
$2\times 3$ matrix of linear forms; an examination of cases shows that $V$ is
Cohen-Macaulay of height two, so is  determinantal. Thus $V$ has the standard
determinantal minimal resolution.
The last part \eqref{netsqiii} follows immediately from Lemma
\ref{netsq0}\eqref{netsq0iii}, and a computation in {\sc Macaulay}.
\end{proof}\par
\begin{lemma}\label{deg2}\cite[Section 3.4-3.6]{Lee} The Hilbert scheme
$\Hilb ^{2,-2}(\mathbb P^3)$
parametrizing curves
$C\subset
\mathbb P^3$ of degree 2, genus -2 (Hilbert polynomial $2t+3$) has two
irreducible components. A general
point of the first parametrizes  a scheme
consisting of two skew lines union a point off the line; this component has
dimension 11. A general point
of the second component parametrizes a planar conic union two points; this
component has dimension 14. \par
Likewise, \cite[Theorem 3.5.1]{Lee} $\Hilb^{2,-1}(\mathbb P^3)$ (Hilbert
polynomial $2t+2$) has the
analogous components parametrizing
two skew lines, or a planar conic union a point. The scheme
$\Hilb^{2,0}(\mathbb P^3)$ (Hilbert polynomial $2t+1$)
has a single component, whose generic points parametrize plane conics.\par
\end{lemma}
The following result mostly concerns certain ideals $I$ for which $I_3$ to
$I_4$ or $I_4$ to $I_5$ is of extremal
growth in the sense of F.H.S. Macaulay. We thank a referee for
the simple argument
for \eqref{Hrestii}. Note that nets $V$ with no linear relation need not
define complete
intersections, and the ideal $(V)$ need not be saturated: thus
\eqref{Hrestiii} below does not
follow from \eqref{Hrestii}.
\begin{lemma}\label{Hrest} Assume for \eqref{Hresti},\eqref{Hrestii} below
that $I$ is a saturated
ideal of
$R=K[w,x,y,z]$.
\begin{enumerate}[i.]
\item\label{Hresti} If
$H(R/I)=(1,4,7,10,13,16,\ldots )$, then $ I_{\le 3}$ defines a twisted cubic
(or specialization not in
the closure of the plane cubics) or a
plane cubic union a point (possibly embedded).  In the former case, $I_2$
lies in $\mathfrak{F}_2$, and generates
$I$; in the latter case
$I_2\in
\mathfrak{{F}_3}$.\par
\par
\item\label{Hrestii}  $H(R/I)$ cannot be any of $(1,4,7,8,10,\ldots
), (1,4,7,b,9,11,\ldots )$, or $ (1,4,7,9,12,\ldots $).
\item\label{Hrestiii}
 If $R/I$ is Artinian Gorenstein of socle degree at least 5, then $R/(I_2)
$ cannot have a
Hilbert function of the form  $H(R/(I_2))= (1,4,7,8,10,\ldots
)$, $H(R/(I_2)) = (1,4,7,b,9,11,\ldots )$, or
$H(R/(I_2))=(1,4,7,9,12,\ldots )$.
\end{enumerate}
\end{lemma}
\begin{proof}Suppose that a saturated ideal $I$ has the Hilbert function
given in case \eqref{Hresti}.  Then 13 to 16 is an extremal growth.  So, by
 the Gotzmann theorem
$I$ defines a scheme
$\Z\subset
\mathbb P^3$, of Hilbert polynomial $3t+1$ so $\Z$
is a degree three curve of genus zero.
The Piene-Schlessinger
Theorem characterizing the components of $\Hilb^{3,0}(\mathbb{P}^3)$
\cite{PS} implies that if
$\Z$ is non-degenerate (not contained in a plane), then $\Z$ is either a
twisted cubic or a specialization, so
$I_2$ is in
$\overline{\mathfrak{F}_2}$, or $\Z$ is the union of a planar cubic and a
(possibly
embedded) spatial point, and then
$I_2$  is in $\mathfrak{F}_3$. If $\Z$ is degenerate, then also $I_2\in
\mathfrak{F}_3$. This completes the proof of \eqref{Hresti}.\par
The three sequences of \eqref{Hrestii}  cannot occur for a saturated ideal
$I$: a saturated
ideal has depth at least one, so $A=R/I$ has a (linear) non-zero divisor,
and the first differences
$\Delta H(R/I)$ must be admissible. But $(1,3,3,1,2,..),
(1,3,3,b-7,9-b,2,..)$ and
$(1,3,3,2,3,..)$ are not $O$-sequences. \par
  In the first case of \eqref{Hrestiii} we have that  $10=8^{(3)}$, so by
Theorem
\ref{MacGo}\eqref{Gotzper} $\Z=\Proj(R/(I_3))$ is a scheme of Hilbert
polynomial $2t+2$ (degree
two and genus -1) and regularity degree no more than 3, the Gotzmann
regularity degree of $2t+2$.
By a classical degree inequality, such a scheme is either reducible, or
degenerate --- contained
in a hyperplane \cite[p. 173]{GH}. Furthermore, by Lemma \ref{deg2} the
Hilbert scheme
$\Hilb^{2,-1}(\mathbb P^3)$ of degree two genus -1 curves has two
irreducible components, one
whose generic point parametrizes two skew lines, the second, whose generic
point parametrizes a planar conic
union a point. For either component, the Hilbert function
$H(R/I_\Z))_2\le 6$ which by Corollary \ref{include} implies $H(R/I)_2\le
6$, contradicting the
assumption. A similar argument handles the second case of \eqref{Hrestiii}:
since $9^{(4)}=11$,
$H_{4,5}=(9,11)$ is maximal growth; by Theorem \ref{MacGo} \eqref{Gotzper}
the scheme
$\Z=\Proj(R/(I_4))$ has Hilbert polynomial
$2t+1$, of Gotzmann regularity two implying
$H(R/I_\Z)_2=5$, and by Corollary \ref{include}, $H(R/I)_2\le 5$, a
contradiction.
For the last case it suffices by Corollary \ref{include} and the Gotzmann
Theorem to know that
any scheme of Hilbert polynomial $3t$ (degree three and genus one) is a planar
cubic or degenerate,
a result of the classification of curves \cite{PS, Lee}.
\end{proof}\par
\subsection{Ideals with $I_2=\langle wx,wy,wz\rangle$.}\label{struc}

 Let $\mathfrak{V}$ denote the vector space $\langle
wx,wy,wz\rangle$. In this section we assume $H=(1,4,7,\ldots ,1)$ and we
consider the subfamily
${\mathfrak C}(H)\subset \PGOR(H)$ parametrizing those algebras $A=R/I$ of
Hilbert function $H$ for
which
 $I_2$ is $\PGL (3)$ isomorphic to
$\mathfrak V$. We first determine when $\mathfrak C(H)$ is nonempty and
give a structure theorem
for such $A$ (Theorem
\ref{V,W}). We then determine the minimal resolution of
$A$ (Theorem
\ref{ResI,J}). We also
 determine the tangent
space to the family $\mathfrak{C}(H)$ (Theorem \ref{sevcomp}). To prove our
results we connect these Artinian algebras with height three Artinian
Gorenstein quotients $R'/J_I$
of
$R'=K[x,y,z]$, where $J_I=I\cap R'$, which are well understood
\cite{BE,D,Klj2,IK}.
\par We recall from
Lemma~\ref{MacD}ff. that, given an ideal $I$ of
$R$, we denote by $I^\perp$ its inverse system, the perpendicular
$R-$submodule to $I$ in the divided power ring $\mathcal
D=K_{DP}[W,X,Y,Z]$, where $R$
acts by contraction.
\begin{theorem}\label{V,W} Let $H=(1,4,7, \ldots )$ of socle degree $j\ge
4$ be a Gorenstein
sequence, and assume that
$I\in
\mathfrak C(H)$ satisfies $I_2=\mathfrak V=
\langle wx,wy,wz\rangle
$.  Let $F\in
\mathcal D_j$ satisfy
$I=\Ann(F)$.
 Let $R'=K[x,y,z]$.
Then,
\renewcommand{\theenumi}{\roman{enumi}}
\begin{enumerate}[i.]
\item\label{V,Wi} The inverse system $(\mathfrak{V})^{\perp}$ of
the ideal
$(\mathfrak{V}),\mathfrak{V}=\langle wx,wy,wz\rangle \subset R$, satisfies
\begin{equation}\label{Vdual}{(\mathfrak{V})^{\perp}}_j=\langle
K_{DP}[X,Y,Z]_j,W^{[j]}\rangle.
\end{equation}
\item\label{V,Wii}     $F \in K_{DP}[W,X,Y,Z]_j$ and satisfies
\begin{equation}\label{eqvspduali}
F=G+\alpha\cdot W^{[j]}, \quad G\in  K_{DP}[X,Y,Z]_j,  \alpha\in K,
\end{equation}
where $G\ne 0, \alpha\ne 0$.

Furthermore,
$I=(J_I,\mathfrak{V},f)$ where
$J_I=I\cap R'$ is the height three Gorenstein ideal $Ann_{R'}(G)$
and $f=w^j-g, g\in K[x,y,z]_j, g\not= 0$.

 The Hilbert function $H(R/I)_i=H(R'/J_I)_i+1$ for $1\le i \le j-1$, so we have
\begin{equation}\label{VIHilb}
H(R/I)=H(R'/J_I)+(0,1,1,\ldots ,1,0)=(1,4,\ldots ,4,1).
\end{equation}
The inverse system
$I^{\perp}$ satisfies $I^{\perp}_j=\langle F\rangle$, $I^{\perp}_{i}=0$ for
$i\ge j+1$, and
\begin{equation}\label{VIdual}
I^{\perp}_i=(R\circ F)_i=\langle (R'\circ G)_i,W^{[i]}\rangle \text{ for }
1\le i\le
j-1.
\end{equation}
\item\label{V,Wiiia} The Gorenstein sequence $H=(1,4,7,\ldots )$ satisfies
$\mathfrak{C}(H)$ is nonempty if and only if $H'=H-(0,1,1,\ldots ,1,0)$ is
a Gorenstein sequence of
height three. (See Corollary \ref{nonempty3}).
\end{enumerate}
\end{theorem}
\begin{proof} We first prove \eqref{V,Wi}. Since
$\mathfrak{V}=(wx,wy,wz)=w\cap (x,y,z)$ we have from the
properties of the Macaulay duality,
\begin{equation*}
(wx,wy,wz)^{\perp}=(w)^{\perp}+(x,y,z)^{\perp}=K_{DP}[X,Y,Z]+K_{DP}[W],
\end{equation*}
which is \eqref{Vdual}.

We now show \eqref{V,Wii}. Since $F$ generates
$(I_j)^\perp$,  $F\in
{(\mathfrak{V})^{\perp}}_j$ can be written
$F=G+\alpha W^{[j]}$ as in \eqref{eqvspduali}.   Since $H(R/I)=(1,4,\ldots
)$, we have
$G\ne 0$ and $\alpha \ne 0$. The inverse system relation \eqref{VIdual} is
immediate, and gives
\begin{equation*}
R\circ F={R'}_{\ge 1}\circ
h+\langle W,W^{[2]},\ldots ,W^{[j-1]},F\rangle,
\end{equation*}
as well as the
Hilbert function equality
\eqref{VIHilb}.  Let $J_I=\Ann(G)\cap K[x,y,z]$: evidently,
$\Ann(G)=(w,J_I)$. Let $h\in I\cap K[x,y,z]$. Then
we have $h\circ F=0$ and $h\circ W^j=0$, implying $h\circ G=0$ so $h\in
J_I$; conversely, if $h\in J_I=
\Ann(G)\cap K[x,y,z]$ then $h\circ G=0, h\circ W^j=0$, implying $h\circ
F=0$, so $h\in I\cap K[x,y,z]$. Thus
$J_I=I\cap K[x,y,z]$, as claimed.  as specified in
\eqref{V,Wii} is immediate, Now,  any form $h$  of degree less than $j$
satisfying $h\cdot F=0$, and $h\ne
(wx,wy,wz)$ must satisfy $h\in K[x,y,z]$ and hence is in $J_1$.
 If   $f=w^j-g$ with
$g\circ G=\alpha$ then we have $f\circ F=0$ and hence $f\in I$.  If
$g= 0$ we would have
$R_1\cdot w^{j-1}\in I$, implying that $w^{j-1}\mod I$ is a socle element
of $A=R/I$, contradicting the
assumption that
$A$ is Artinian Gorenstein of socle degree $j$. Thus, we have $f=w^j-g$
with $g\ne 0$. Since the lowest-degree
third syzygy of
$I$ are those in degree four arising from
$\mathfrak{V}$, the symmetry of the minimal resolution implies that $I$ has
no generators (first syzygies) in
degrees greater than $j$. Thus
the ideal
$I\in
\mathfrak{F}$ is minimally generated as
$I=(J_I,\mathfrak{V},f)$.  as claimed, and completes
the proof of \eqref{V,Wii}.

To show \eqref{V,Wiiia}, note that if $I\in
\mathfrak {C}(H)$ then $H'$ from
\eqref{V,Wiiia} satisfies
$H'=H(R/J_I)=H(R'/(I\cap R'))$ with $I\cap R'$ a Gorenstein ideal in $R'$, so
$H'$ is a Gorenstein sequence. Conversely if $H'=H-(0,1,1,\ldots ,1,0)$ is
a Gorenstein sequence
then  take $J$ to be any Gorenstein ideal in $R'$ of Hilbert function $H'$ and
let $J = \Ann _{R'}(G)$.  Let $F = G+W^j$.  Then
$\Ann (F) = I= (J,w^j-g,wx,wy,wz)$
where $g\in R'_j$ but $g\notin J$: the ideal $I$ is a Gorenstein ideal of
height four.
  Then  we have $I\in
\mathfrak {C}(H)$.

Thus, $\mathfrak {C}(H)$ is nonempty if and only if
$H'=H-(0,1,1,\ldots , 1,0)=(1,3,\ldots ,3,1)$ is a
Gorenstein sequence of height three.  This completes the proof.
\end{proof}
\par
The minimal resolution of $R/I$ can be constructed from the minimal
resolution of $J_I$.  We construct a putative complex in Definition
\ref{resolveI};
we prove that it is an exact complex in Theorem~\ref{ResI,J}.  The
construction relies on Theorem \ref{V,W}\eqref{V,Wii}.\par
 Suppose that $I\subset R$ defines a Artinian Gorenstein quotient $A=R/I$,
that $I_1=0$ and
$I_2=\mathfrak{V}$, and that
$I=(\mathfrak{V},J_I,g-w^j)$ with $g\in R'_j$ satisfying $g\ne 0$, and
$J=J_I=I\cap R'=K[x,y,z]$ defining a
Artinian Gorenstein quotient $A'=R/J'$ of $R'$. Let the minimal resolution of
$R/J$ be (here $m=2n+1$ is odd)
\begin{equation}\label{Jres}\mathbb{J} :\quad 0\to
R\xrightarrow{\alpha^t}R^{m}\xrightarrow{\phi}R^{m}\xrightarrow{\alpha} R\
\to R/J\to 0,
\end{equation}
where $\phi$ is an $m\times m$ alternating matrix with homogeneous
entries, and $\alpha=\left[ J\right]$
denotes the $1\times m$ row vector with entries the homogeneous
generators of $J$ that are the
Pfaffians of
$\phi$, according to the Buchsbaum-Eisenbud structure theorem for height
three Gorenstein ideals (since $J$ is
homogeneous, $\mathbb{J}$ may be chosen homogeneous: see
\cite{BE,D}). Denote by
$\mathbb{K}$ the Koszul complex resolving $R/(x,y,z)$ (so $\mathbb K_0=
\mathbb K_3=R$):
\begin{equation}
\mathbb{K} :\quad 0\to R \xrightarrow{\delta_3} R^3 \xrightarrow{\delta_2}
R^3 \xrightarrow{\delta_1}R\to
R/(x,y,z) \to 0,
\end{equation}
where $\delta_1=[x,y,z], \delta_2=\left( \begin{smallmatrix}
y&z&0\\
-x&0&z\\
0&-x&-y
\end{smallmatrix}\right)$, and  $\delta_3=\delta_1^t$. We will let
$\mathbb{T}:\mathbb{K}\to   \mathbb{J}$  be a
map of complexes induced by multiplication by $g$ on $R$. By degree
considerations, we see that
$\deg T_3=0$, so $T_3$ is multiplication by $\gamma\in K$. \par
So we have $T_1\circ \delta_2=\phi\circ T_2$, also $T_2\circ
\delta_3=[J]^t$, and
\begin{equation*}T_2\circ  \left[ \begin{smallmatrix}z&\\ y&\\ x&
\end{smallmatrix} \right] = \gamma \circ
[J]^t.
\end{equation*}
\begin{definition}\label{resolveI} Given $I,J,\mathbb J, \mathbb K$ as above,
 we define the following complex,
\begin{equation}\label{complexF}\mathbb{F}:\quad 0\to
R\xrightarrow{F_4}R^{m+4}\xrightarrow{F_3}R^{2m+6}\xrightarrow{F_2}R^{m+4}
\xrightarrow{F_1}R\to R/I\to 0,
\end{equation}
where $F_1=(wx,wy,wz,\alpha,w^j-g)$, and $F_2$ satisfies
\begin{equation}
F_2\quad=\qquad \left(
\begin{array}{c|cccc}
&3&m&m&3\\ \hline \\
3&\delta_2&0&\frac{-1}{\gamma}E\circ T_2^t&w^{j-1}I_{3\times 3}\\
m&0&\phi &wI_{m\times m}& T_1\\
1& 0& 0&0&-x\ -y\ -z
\end{array}\right),
\end{equation}
where $E=\left[
\begin{smallmatrix}
0&0&1\\
0&1&0\\
1&0&0
\end{smallmatrix}\right]$. The map $F_3$ satisfies
\begin{equation}
F_3\quad=\qquad \left(
\begin{array}{c|ccc}
&3&m&3\\ \hline \\
3&w^{j-1}I_{3\times 3}&ET_1^t&-z\ -y\ -x\\
m&T_2&-{\gamma}wI_{m\times m}&0\\
m&0&\gamma\phi &0\\
3& -\delta_2& 0&0
\end{array}\right),
\end{equation}
 and $F_4=(wz,wy,wx,\alpha,w^j-g)^t$.
\end{definition}

\begin{theorem}\label{ResI,J}
Let I be a homogenous height four Gorenstein ideal in $R = K[w,x,y,z]$ with
socle degree j and with  $I_2 = (wx,wy,wz)$. Then the complex $\mathbb{F}$ of
\eqref{complexF} in Definition \ref{resolveI} is exact and is the minimal
resolution of $R/I$.
\end{theorem}
\begin{proof}
  We first show that
$\mathbb{F}$ is a complex. By \eqref{V,Wii} of the structure theorem, we
see that
$I$ is minimally generated by $J = I\cap K[x,y,z], wx,wy,wz , g-w^j $ where
$g\in K[x,y,z]$.  So, $g\notin J$.  Suppose that $\gamma=0$. Then
$T_2\circ \delta_3=0$, hence we would have
$T_2=T'\circ\delta_2$ for some
$T'$. Then
\begin{align*}
T_1\circ\delta_2&=\phi\circ T_2=\phi\circ T'\circ\delta_2;\\
\mbox{so              }\qquad (T_1-\phi\circ T')\circ\delta_2&=0;
\mbox{and}\\
T_1-\phi\circ T'&=\beta [x,y,z], \beta\in K,\mbox{  and}\\
\alpha\circ T_1&=\alpha\circ \beta [x,y,z], \mbox{ and}\\
-g[x,y,z]&=\alpha\beta [x,y,z].
\end{align*}
This implies $g\in J$ contradicting $g\notin J$.  So, we get
 $\gamma\ne 0$.
\par We get $F_1\circ F_2=0$ and $F_3\circ F_4=0$ from the following three
identities. First, from the exact
sequence $\mathbb J$ of \eqref{Jres} we have
\begin{equation} \phi\alpha=\alpha^t\phi=0.
\end{equation}
Second, from
\begin{align} [x,y,z]\left( (\frac{-1}{\gamma})ET_2^t\right)
&=\frac{-1}{\gamma}[z y x]T_2^t\notag\\
&=\frac{-1}{\gamma}\left[ T_2\left[
\begin{smallmatrix}
z&\\
y&\\
x&
\end{smallmatrix}
\right]
\right]^t=\frac{-1}{\gamma}\gamma (\alpha^t)^t\notag\\
&=-\alpha\notag\\
\mbox{ we have }\qquad [x,y,z]\left[ \frac{-1}{\gamma}ET_2^t\right] &
=-\alpha .
\end{align}
Third, we have
\begin{equation} T_1J=-g[x,y,z] .
\end{equation}
To see that $F_2\circ F_3=0$ we just need to check that
\begin{align*}
\phi T_2-T_1\delta_2&=0\\
\mbox{and}\qquad \delta_2ET_1^t-\frac{1}{\gamma}ET_2^t(\gamma \phi)&=0.
\end{align*}
The first of these follows from the map of complexes
$\mathbb{T}:\mathbb{K}\to\mathbb{F}$. For the second\
we have
\begin{align*}\delta_2ET_1^t-ET_2^t\phi&=\delta_2ET_1^t+ET_2^t\phi^t\\
&=\delta_2ET_1^t+E(\phi T_2)^t\\
&=\delta_2ET_1^t+E(T_1\delta_2)^t\\
&=\delta_2ET_1^t+E(\delta_2^t)T_1^t\\
&=\left( \delta_2E+E\delta_2^t\right) T_1^t = 0,\\
\mbox{since}\qquad &\delta_2E+E\delta_2^t=0.
\end{align*}
So we get $F_2F_3=0$. Thus, $\mathbb{F}$ is a complex.\par To see that
the complex $\mathbb{F}$ is exact, we
use the exactness criterion \cite{BE1}\cite[Theorem 20.9]{Ei}. It suffices
to show that $\sqrt{I_{m+3}(F_2)}$
and
 $\sqrt{I_{m+3}(F_3)}$ have depth at least three, where $I_{m+3}(F_2)$
denotes the Fitting ideal generated by
the
$(m+3)\times (m+3)$ minors of $F_2$.
We write $F_2$ as
\begin{equation}
F_2\quad=\qquad \left(
\begin{array}{c|cccc}
&3&m&m&3\\
 \hline \\
3&\begin{smallmatrix}y&z&0\\-x&0&z\\0&-x&-y\end{smallmatrix}&0&
\begin{smallmatrix}\tau _{11}&\ldots
&\tau _{1m}\\&\ldots&\\
\tau_{31}&\ldots&\tau_{3m}\end{smallmatrix}&
\begin{smallmatrix}w^{j-1}&0&0\\0&w^{j-1}&0\\0&0&w^{j-1}
\end{smallmatrix}\\
m&0&\phi
&\begin{smallmatrix}w\ &0\ &0\ \\0\ &w\ &0\ \\0\ &0\ &w\ \end{smallmatrix}
&T_1\\ 1& 0& 0&0&-x\ -y\ -z
\end{array}
\right),
\end{equation}
where $x\tau_{1i}+y\tau_{2i}+z\tau_{3i}=-\alpha_i$, and $J=(\alpha_1,\ldots
,\alpha_m)$. Consider the
minor
$M_i$ of $F_2$ having all rows except the $(3+i)$-th row, and having the
columns $1,2,4,\ldots
,3+i-1,3+i+1,m+3,m+3+i,2m+4$. This is the minor
\begin{equation}
M_i\quad=\quad\qquad
\begin{array}{|cccc|}
\begin{smallmatrix}y&z\\-x&0\\0&-x\end{smallmatrix}&0&
\begin{smallmatrix}t_{1i}\\t_{2i}\\t_{3i}\end{smallmatrix}&\begin{smallmatrix}w^
{j-1}\\0\\0\end{smallmatrix}\\
0&\phi_i &0&\begin{smallmatrix}\ast \\\ast \\\ast \end{smallmatrix}\\
0& 0&0&-x
\end{array},
\end{equation}
and it equals
\begin{align*}&\pm
xa_i^2\begin{vmatrix}y&-z&t_{1i}\\-x&0&t_{2i}\\0&-x&t_{3i}\end{vmatrix}\\
&=\pm xa_i^2x(x\tau_{1i}+y\tau_{2i}+z\tau_{3i}\\
&=\pm x^3a_i^2.
\end{align*}
Thus $xa_i\in \sqrt{I(F_2)}$. Similarly, $ya_i,za_i\in \sqrt{I(F_2)}$. Thus
$mJ\subset\sqrt{I(F_2)}$.
Finally, looking at the last $m+3$ rows and the columns $1,2,m+4,\ldots
,2m+4$, we get $\pm x^3w^m$ in
$I(F_2)$. So $wx\in\sqrt{I(F_2)}$, as well as $wy,wz$, by similar
computations. Thus $\sqrt{I(F_2)}\supset
(J,wx,wy,wz)$. Similarly $\sqrt{I(F_3)}\supset (J,wx,wy,wz)$. So these
Fitting ideals have depth at least
three, and the complex $\mathbb F$ is exact.
 This completes the proof.
\end{proof}

\begin{remark} The above resolution in Theorem \ref{ResI,J} is
similar to but different from the
minimal  resolution obtained by A. Kustin and M. Miller in \cite{KuMi2}.
They consider ideals of the form
$ (f, g, h, wJ)$ where  $(f,g,h)$ is a regular sequence and $J$ is height
three Gorenstein. It
turns out that it is not a specialization  of their resolution. One
reason for the  resemblance is
that $(wx,wy,wz)$ has three Koszul type relations even though they are not
a regular sequence.
\end{remark}\par
If $H(R/I)=(1,4,7,h,7,4,1)$, recall that
$\mathfrak {C}(H)\subset
\PGOR(H)$ denotes the subfamily parametrizing ideals $I$ such that
$I_2\cong \mathfrak{V}=\langle
wx,wy,wz\rangle$, up to a coordinate change.  We denote by
$\nu_i(J)$ the number of degree-$i$ generators of $J$.  We will later show
that any
Gorenstein sequence $H=(1,4,7,\ldots )$ satisfies $\mathfrak {C}(H)$
nonempty (Theorem \ref{nonempty2}).
For $I\in \PGOR(H)$  we denote by $\mathcal T_I$ the tangent space
to the
affine cone over $\PGOR(H)$ at the point corresponding to $A=R/I$. Recall that
$H'=H-(0,1,1,\ldots ,1,0)$. We denote by
$\mathcal T_{J_I}$ the tangent space to the affine cone over $
\PGOR(H'), H'=H(R/J_I) $ from \eqref{VIHilb}, at the point corresponding to
$A'=R'/J_I,$ where $ J_I=I\cap
K[x,y,z]$.
\begin{theorem}\label{sevcomp}Let $H=(1,4,7,\ldots)$ of socle degree $j\ge
5$. In
\eqref{sevcompi},\eqref{sevcompii}, \eqref{sevcompiii} we let
$A=R/I\in
\mathfrak {C}(H)$, and we let
$J_I=I\cap K[x,y,z]$.
\begin{enumerate}[i.]
\item\label{sevcompi}
The dimension of $\mathfrak {C}(H)\subset \PGOR(H)$ satisfies
\begin{equation}\label{sevcompieq}
\dim ({\mathfrak {C}(H)})=7+\dim\PGOR(H').
\end{equation}
\item\label{sevcompii} The dimension of the
tangent space
$\mathcal T_I$ to the affine cone over
$\PGOR(H)$ at the point determined by $A=R/I\in\PGOR(H)$ satisfies,
\begin{equation}\label{tangsp0}
\dim_K \mathcal{T}_I=7+\dim_K \mathcal{T}_{J_I}+\nu_{j-1}(J_I).
\end{equation}
\item\label{sevcompiii}  The GA algebra $A\in \mathfrak{E}(H)$ is a smooth
point of $\PGOR(H)$ if and only if
$\nu_{j-1}(J_I)=0$.
\item\label{sevcompiv} The subscheme $
\mathfrak {C}(H)$ of $\PGOR(H)$ is irreducible.
\item\label{sevcompv} When $j=6$ and $H=H_h=(1,4,7,h,7,4,1), 7\le h\le 11$
we have
\begin{equation}\label{dimE}
\dim ({\mathfrak {C}(H)})= 34-\binom{h^\vee +1}{2}, \quad h^\vee = 11-h.
\end{equation}
When also,
$8\le h\le 11$, $\mathfrak {C}(H)$ is generically smooth.
\end{enumerate}
\end{theorem}
\begin{proof} The proof of \eqref{sevcompi} is immediate from the structure
Theorem \ref{V,W}\eqref{V,Wii}: the
choice of $\mathfrak{V}$ involves that of $w$ and the vector space $\langle
x,y,z\rangle$, so 6 dimensions, and
that of the $F=w^j+G$ involves one parameter, given $\langle G\rangle$,
which determines $J_I$. \par
We now show \eqref{sevcompii}. Let $A=R/I\in \mathfrak {C}(H)$. We recall
from \cite[Theorem 3.9]{IK} that for a
GA quotient $A=R/I$, we have
$\dim_K \mathcal T_I=\dim_K R_j/(I^2)_j=H(R/I^2)_j$.
We have
\begin{align*}
 (I^2)_j&= I_2\cdot I_{j-2}\oplus (J^2)_j\\
      &=  \left( wR'_1\cdot
\left((w^{j-3}R'_1\oplus w^{j-4}R'_2\oplus\cdots \oplus wR'_{j-3})\oplus
J_{j-2}\right)\right)\oplus (J^2)_j\\
&=\left( w^{j-2}R'_2\oplus w^{j-3}R'_3\oplus \cdots \oplus
w^2R'_{j-2}\right)\oplus wR'_1J_{j-2}\oplus (J^2)_j.
\end{align*}
Hence we have
\begin{align*}
 R_j/(I^2)_j&\cong
w^j\oplus w^{j-1}R_1'\oplus w\left( R'_{j-1}/R'_1J_{j-2}\right)\oplus
R'_j/(J_j)^2, \mbox{   and }
\\ \dim_K R_j/(I^2)_j&=1+3+H'_{j-1}+\nu_{j-1}(J)+\dim_K R'_j/(J_j)^2\\
&=7+\dim_K \mathcal{T}_{J_I}+\nu_{j-1}(J_I).
\end{align*}
We now show \eqref{sevcompiii}. We
use J.-O. Kleppe's result that in codimension 3,
$\PGOR(H')$ is smooth \cite{Klj2}. It follows that for the Gorenstein ideal
$J_I\subset R'=K[x,y,z]$, of socle
degree
$j$, of
 Hilbert function $H(R'/J)=H'$
the dimension of the tangent space $\mathcal T_{J_I}$ to the affine cone
over $\PGOR(H')$ at $J_I$ satisfies
\begin{equation*}
\dim_K T_{J_I}=\dim (\PGOR(H'))+1.
\end{equation*}
This, together with \eqref{sevcompi},\eqref{sevcompii} shows that
$\nu_{j-1}(J_I)=0$ implies $\dim_K \mathcal
T_I=\dim \mathfrak{E}(H)+1$, hence that $\mathfrak{E}(H)$ and $\PGOR(H)$
are smooth at such points, which is
\eqref{sevcompiii}.\par
We now show \eqref{sevcompiv}. We first show that $\mathfrak {C}(H)$ is
irreducible. The scheme $\PGOR(H')$ is
irreducible by
\cite{D} (or by its smoothness
\cite{Klj2}, discovered later). The scheme $\mathfrak {C}(H)$, is fibred
over the family of nets isomorphic
to
$\mathfrak{V}$ by $\PGOR(H')$, then by an open in $\mathbb P^1$ (to choose
$F$ given $G$), so it is
irreducible.

\par
 We now show \eqref{sevcompv}. The dimension formula \eqref{dimE} results
immediately from \eqref{sevcompi}
and the known dimension of $\PGOR(H')$ (see \cite[Theorem
4.1B]{IK},\cite{Klj2}). From the latter source, we have
that the codimension of $\PGOR(H')\subset \mathbb P^{27},
H'=(1,3,6,h-1,6,3,1)$ is $\binom{h^{\vee}+1}{2}$ where
$h^\vee=10-(h-1)$. When also $8\le h\le 11$, we have $\Delta^3(H')_5=0$; it
follows simply from \cite{D} (or see
\cite[Theorem 5.25]{IK}) that the generic GA quotient  $R'/J$ having
Hilbert function $H'$ satisfies $\nu_5(J)=0$.
This completes the proof of
\eqref{sevcompv} and of the Theorem.
\end{proof}\medskip
\par
\subsection{Mysterious Gorenstein algebras with $I_2=\langle
w^2,wx,wy\rangle$}\label{w2sec}

 Let $\mathfrak{W}$ denote the vector space $\langle
w^2,wx,wy\rangle$. In this section we assume $H=(1,4,7,\ldots ,1)$ and
study  graded Artinian
Gorenstein algebras
$A=R/I, R=K[w,x,y,z]$, such that
\begin{equation}\label{ebase}
A\in \mathfrak{E}_{sp}(H):\
I_2=\W.
\end{equation}
We will show that their Hilbert functions are closely related
to those of a
Gorenstein ideal in three variables (Lemmas \ref{HFw1},\ref{HFw2}). From
these results we can characterize the Hilbert functions $H$ for which
$\mathfrak E_{sp}(H)$ is
nonempty (Theorem
\ref{hfw2}): these are the same as found in the previous
section for Gorenstein algebras
$A\in \mathfrak C(H)$: those with $I_2\cong \langle wx,wy,wz\rangle$. However
 it is
an open question whether the Zariski closure $\overline{\mathfrak C (H)}$
contains
$\mathfrak E_{sp}(H)$, and it is this uncertainty that requires us to
consider $\mathfrak E_{sp}$
in detail.
\par
The ideal $(\W)$ generated by $\W$ satisfies $(\W)=(w^2,x,y)\cap
(w)$.  The inverse system $\W^\perp\subset \mathcal D$ satisfies
\begin{align}\label{eW}\W^{\perp}=\left( (w^2,x,y)\cap
(w)\right) ^{\perp}=&(w^2,x,y)^{\perp}+(w)^{\perp}\\ =&\notag
K_{DP}[Z]+W\cdot K_{DP}[Z] +K_{DP}[X,Y,Z],
\end{align}
Thus we have for the degree-$j$ component
\begin{equation*}
\{\W^\perp\}_j=K_{DP}[X,Y,Z]_j+\langle
WZ^{[j-1]}, Z^{[j]}\rangle .
\end{equation*}
\begin{lemma}   Let $I$ satisfy \eqref{ebase}, and let
$F\in \mathcal{R} =K_{DP}[W,X,Y,Z]_j$  be a generator
of its inverse system. Then $F$ may be written uniquely
\begin{equation}\label{esum}
F=G+WZ^{[j-1]}, \, G\in K_{DP}[X,Y,Z],
\end{equation}
in the sense that the decomposition depends only on $I$, and the choice of
generators $w,x,y,z$ of $R$. Further, after a linear change of basis in $R$, we
may suppose that
$G$ in
\eqref{esum} has no monomial term in $Z^{[j]}$.
\end{lemma}
\begin{proof}   Since $w^2, wx, wy$ are all in $I$, by \eqref{eW} the
generator $F$ of $I^{\perp}$
can be written in the form $F=G+\lambda WZ^{[j-1]}, G\in K_{DP}[X,Y,Z]$.
Evidently, $\lambda \ne 0$, since otherwise $H(A)=(1,3,\ldots
)$; so we may choose $\lambda =1$. The decomposition of \eqref{esum} is
certainly unique, given $I$, and the choice of $x,y,z,w$.
 A linear change of basis $w\to w, x \to x, y\to
y, z\to z+\beta w$ in
$R$, and the contragradient change of basis $W\to W-\beta Z,
X\to X, Y\to Y, Z\to Z$ in
$\mathcal R$ eliminates any
monomial term in $Z^{[j]}$ from $G$.
\end{proof}\par
We denote by $R'$ the polynomial ring $R'=K[x,y,z]$.
\begin{lemma}\label{alpha} Let $I$ be an ideal satisfying
\eqref{ebase}, let
$F=G+WZ^{[j-1]}$ be a generator of its inverse system as in
\eqref{esum}, and let $J=\Ann(G), J'=\Ann(G)\cap R'$; then $J=(w,J')=J'R$.
Let
$\alpha (J)$ be the integer
\begin{equation}\label{ealpha}
\alpha (J)=\min\{ \alpha\ge 1 \mid J'_\alpha \nsubseteq (x,y)\}=\min\{ \alpha\ge
1 \mid J_\alpha \nsubseteq (x,y,w)\}.
\end{equation}
Then $\alpha(J)=\min\{ i\mid Z^{[i]}\notin R'_{j-i}\circ G\}$, and
we have $2\le \alpha(J)\le  j$.
\end{lemma}
\begin{proof} The first statement follows from $(x,y)^\perp\cap
\mathcal{R}'=K_{DP}[Z]$. The lower bound on $\alpha$ follows
from the assumption of \eqref{ebase}, which implies that
$H(R/J')=(1,3,\ldots )$, so $2\le \alpha$. The upper bound on
$\alpha$ follows from the fact that   $z^j \in J' = \Ann(G)$.
\end{proof}
\begin{definition}\label{halpha}Let $I$ satisfy \eqref{ebase},
let
$F=G+WZ^{[j-1]}$ be a generator of its inverse system, as in
\eqref{esum}, and let $\alpha=\alpha(J)$ as in \eqref{ealpha}.
We define a sequence
\begin{equation}\label{eH}
 H_\alpha=
\begin{cases}(0,1,1,\ldots ,1,
2=h_\alpha,2,\ldots ,2=h_{j-\alpha},1,\ldots ,1,0=h_j)
\text{ if $\alpha\le
j/2$,}\\
(0,1,1,\ldots ,1=h_{j-\alpha},0,\ldots ,0,1=h_\alpha,1,\ldots
1,0=h_j)
 \text{ if $\alpha >j/2$ and $j\ne 2\alpha-1$}\\
(0,1,1,\ldots ,1,0=h_j) \text{ if $j=2\alpha -1$.}
\end{cases}
\end{equation}
We let $H_0=(0,1,1,\ldots , 1,0=h_j)$.
\end{definition}
Note that $H_\alpha$ takes values only 0,1, and 2. When
$\alpha\le j/2$, there are
$j+1-2\alpha
$ 2's in the middle of the sequence $H_\alpha$; when
$\alpha>j/2$ there are $2\alpha+1-j$ 0's
in the middle of
$H_\alpha$. When $\alpha \le j/2$ the middle run of 2's is bordered on the left
by 0 in degree zero, followed by $\alpha-1$ 1's. When $\alpha >j/2$ the
middle run of 0's is bordered on the left by 0 in degree zero followed by 
$j-\alpha$ 1's.\par
\begin{definition}
 We denote by $M$ the $R$-submodule of $\mathcal D$ generated by
$WZ^{[j-1]}$, whose degree-$i$ component satisfies $M_i=\langle
Z^{[i]},W\cdot Z^{[i-1]}\rangle$ for $1\le i<j$. Given $F,G$ as in
\eqref{esum} we define two $R$-modules
\begin{align}
B=&R\circ
\langle F,WZ^{[j-1]}\rangle /R\circ G \notag\\ C=&R\circ
\langle F,WZ^{[j-1]}\rangle /R\circ F.\label{eBC}
\end{align}
We denote by $H^\vee (B)$ the
\emph{dual} sequence $H^\vee (B)_i=H(B)_{j-i}$,
and likewise $H^\vee (C)_i=H(C)_{j-i}$.
\end{definition}\noindent
 Evidently we have for $F,G$ as in \eqref{esum}
\begin{equation}\label{eFrel} I\cap J=\Ann \langle F,G\rangle =\Ann
\langle F,WZ^{[j-1]}\rangle =\Ann \langle G,WZ^{[j-1]}\rangle .
\end{equation}
Our convention will be to specify Hilbert functions of $R$-submodules of
$\mathcal D$ (or
of $\mathcal R$) as subobjects: thus $H(R\circ \{Z^{[2]},WZ\})=H(\langle
1;Z,W;Z^{[2]},WZ\rangle)=(1,2,2)$. However, the Hilbert functions $H(B)$,
and $H(C)$ are as
$R$-modules: thus, when $F=X^{[2]}\cdot Z^{[2]}+WZ^{[3]}$, the module $B$ from
\eqref{eBC} satisfies, after taking representatives for the quotient,
$
B\cong\langle WZ^{[3]};Z^{[3]},W\cdot Z^{[2]}; WZ;W\rangle$
 so $H(B)=(1,2,1,1),
$
and the dual sequence $H'(B)=(0,1,1,2,1)$.
\begin{lemma} We have
\begin{align}
H(R/(I\cap J))=&H(R'/J')+H^\vee (B)\notag\\
=&H(R/I)+H^\vee (C).\label{ehilbrel}
\end{align}
  The $R$-modules $B$ and $C$
each have a single generator, the class of $WZ^{[j-1]}$.
\end{lemma}
\begin{proof} \eqref{ehilbrel} is immediate from \eqref{eFrel}, and the
definition of $H(B),H(C)$. The last statement is immediate from the
definition of $B,C$.
\end{proof}
\begin{lemma}\label{HFw1} Let $I$ be an ideal satisfying
\eqref{ebase}, and let
 $F=G+WZ^{[j-1]}$ be a decomposition as in \eqref{esum} of
the generator $F$ of the inverse system $I^{\perp}$.  Let $J = \Ann (G)$
and  $\alpha =\alpha (J)$ as in
\eqref{ealpha}. Then we have
\begin{enumerate}[i.]
\item \label{HFw1i} $I\cap J=\Ann\langle G,WZ^{[j-1]}\rangle$, and
$(I\cap J)^{\perp}=\langle R'\circ G,M\rangle
=(J')^{\perp}+M=I^{\perp}+M$.
\item\label{HFw1ii} $H(B)=(1,2,\ldots ,2_{j-\alpha},1,\ldots
,1,0)$,
$H(C)=(1,1,\dots 1_c)$, with $c=\alpha$ or $c=j-\alpha$. The
case $c=j-\alpha$ can occur only if $\alpha \ge j/2$.
\item\label{HFw1iii} When $\ c=\alpha$, we have
$H(R/I)-H(R'/J')=H_\alpha$; when $c=j-\alpha$ we have
$H(R/I)-H(R'/J')=H_0=(0,1,1,\ldots ,1,0)$.
\end{enumerate}
\end{lemma}
\begin{proof} Since  $I=\Ann(F)=\Ann(G+WZ^{[j-1]})$ and
$J=\Ann(G)$, we have
\begin{equation*}
I\cap J=\Ann\langle F,G\rangle =\Ann\langle
G, WZ^{[j-1]})=\Ann\langle F, WZ^{[j-1]}\rangle.
\end{equation*}
This proves \eqref{HFw1i}. To show \eqref{HFw1ii} we consider
the two $R$-modules $B,C$ defined above. Evidently we have
 $H(B)_i\le 2$, whence by the Macaulay inequalities
$H(B)=(1,2,\ldots ,2_a,1,\ldots ,1_b,0)$, with invariants the
length $a-1$ of the sequence of $2's$, and the length $b-a$ of
the sequence of $1's$. Since $Z^{[i]}\in R\circ F$ for $1\le i\le
j-1$, we have that the Hilbert function $H(C)$ satisfies
$H(C)_i
\le 1$, hence,
$H(C)=(1,1,\ldots 1_c,0)$, with sole invariant the length
$c+1$ of the sequence of $1's$. Now
\begin{align*} H(R/(I\cap
J))_i-H(R/J)_i=2 & \Leftrightarrow M_i\oplus (R'\circ
G)_i=R_{j-i}\circ \langle G,WZ^{[j-1]}\rangle \\
&\Leftrightarrow Z^{[i]}\notin (R'\circ G)_i\\
&\Leftrightarrow i\ge \alpha (J).
\end{align*}
Otherwise, for $1\le i<\alpha(J), H(R/(I\cap
J))_i-H(R/J)_i=1$, since for such $i$ we have
\begin{equation*}
WZ^{i-1}\in
R_{j-i}\circ
\langle G,WZ^{[j-1]}\rangle \text{ but }
WZ^{[j-1]}\notin R'_{j-i}\circ G,
\end{equation*}
 and for $i=0$ the difference
is
$0$.  Hence, taking into account that $H^\vee
(B)=H(R/(I\cap J))-H(R'/J')$, we have $a=j-\alpha (J)$ and
$b=j-1$. Since both
$H(R/I)$ and
$H(R'/J')$ are symmetric about $j/2$, so is their difference
\begin{equation}\label{ediff}
H(R/I)-H(R'/J')=H^\vee(B)-H^\vee (C).
\end{equation}
This difference can be symmetric only if $c=\alpha$ or
$c=j-\alpha$.  \par
Suppose now that $c=j-\alpha$, and
$\alpha < j/2$. We will show that $H(R/I)\alpha
=H(R/J)_\alpha +2$. By definition of $\alpha$,
$J'_\alpha
$ has a generator of the form
$z^\alpha-g,g\in (x,y)R'$; it follows that $z^{j-\alpha}-g'\in
J', g'=z^{j-2\alpha}g\in (x,y)R'$. Consider the subset
\begin{equation*} \left( (x,y)\cdot R'\right)\circ G = \left(
(x,y)\cdot R'\right)\circ F.
\end{equation*}
Note that $z^{j-\alpha}\circ G\in (x,y)R'\circ G$. However,
$z^{j-\alpha}\circ  F$ has a term $WZ^{[\alpha -1]}$, and
$wz^{j-\alpha -1}\circ F=Z^{[\alpha ]}$. By Lemma \ref{alpha}
$Z^{[\alpha ]} \notin R'\circ G$, it follows that $\dim
R_{j-\alpha}\circ F=\dim R'_{j-\alpha}\circ G +2$, as claimed.
This implies that $H(R/I)=H(R/J)+H_\alpha$ is the only
possibility when $\alpha < j/2$.
\par The statement \eqref{HFw1iii} is immediate from
 \eqref{HFw1ii} and \eqref{ediff}.
\end{proof}
\begin{remark} Note that, {\it given the Hilbert function}
$H'=H(R/J)$ the condition $\alpha (J)\ge  \alpha_0$ is a
closed condition on the familiy $\PGOR (H')$. That is, it is
rarer to have higher values of $\alpha(J)$. However, the situation
is quite different if the Hilbert function is allowed to change, for
example if a
term $\lambda Z^{[j]}$ is added to the dual generator $G$ of $J$: see Lemma
\ref{lambda}, where the effect of such a change is described.
\end{remark}
\begin{lemma}\label{HFw2}Let $I$ be an ideal satisfying
\eqref{ebase}, and suppose that
 $F=G+WZ^{[j-1]}$ be a decomposition as in \eqref{esum} of
a generator $F$ of the inverse system $I^\perp$. Let $\alpha
=\alpha (J), J=\Ann(G)$  be the integer of
\eqref{ealpha}. Then we have
\begin{enumerate}[i.]
\item \label{HFw2i} $H(R/I)$ satisfies either
$H(R/I)=H(R'/J')+H_\alpha$ or
$H(R/I)=H(R'/J)+H_0$; the second possibility may occur only if
$\alpha\ge j/2$.
\item \label{HFw2ii} If $H(R/I)=H(R/J')+H_\alpha$,
then \begin{align*}
H(R/(I\cap J))=&H(R/I)+(0,0,\ldots
,0,1_{j+1-\alpha},1,\ldots ,1_j)
\text{ and }\\(I\cap J)^{\perp}=&I^{\perp}\oplus \langle WZ^{[j-\alpha
]},\ldots
,WZ^{[j-1]}\rangle  \\
\text { also } H(R/(I\cap
J))=&H(R'/J')+(0,1,\ldots ,1,2_\alpha,2,\ldots 2_{j-1},1_j), \text {
and } \\(I\cap J)^{\perp}=&(J')^{\perp}\oplus
\langle W, WZ,\ldots ,WZ^{[j-1]};Z^{[\alpha ]},Z^{[\alpha+1]},\ldots
,Z^{[j-1]}\rangle .
\end{align*}
\item \label{HFw2iii} If $H(R/I)=H(R/J')+H_0$, then $H(R/(I\cap
J))$ and  $H(R'/J')$ are related as above, but
\begin{equation*}
H(R/(I\cap
J))=H(R/I)+(0,0,\ldots ,0,1_\alpha,1,\ldots ,1_j).
\end{equation*}
\begin{proof} The Lemma is an immediate consequence of Lemma
\ref{HFw1} and \eqref{ediff}.
\end{proof}

\end{enumerate}
\end{lemma}
Recall that a Gorenstein sequence $H$ of height $3$ is a
non-negative sequence of integers $H=(1,3,\ldots ,
1=h_j,0,\ldots )$, symmetric about $j/2$, that occurs as the
Hilbert function of a graded Artinian Gorenstein algebra
$A\cong K[x_1,\ldots ,x_r]/I$.  Recall that
then $
 (\Delta H)_i=H_i-H_{i-1}$.\par\noindent
\begin{theorem}\label{hfw2} Let $I$ be an ideal satisfying
\eqref{ebase}. Then $H=H(R/I)$ satisfies
\begin{enumerate}[i.]
\item\label{hfw2i} $\Delta
H_{\le j/2}$ is an $O$-sequence.
\item\label{hfw2ii} $H=H'+H_0=H'+(0,1,1,\ldots ,1,0)$ for some
 Gorenstein sequence $H'$ of height three.
\end{enumerate}
\end{theorem}
Warning:
the $H'$ of \eqref{hfw2ii} above is \emph{not} in general equal to
 $H(R'/J')$, except when $c=j-\alpha$.
\begin{proof} By Lemma \ref{HFw1}\eqref{HFw1ii} we have $c=\alpha $ or
$c=j-\alpha$. The result of the Theorem is obvious in the case
$c=j-\alpha$, since then by Lemma \ref{HFw1}\eqref{HFw1iii}
$H(R/I)=H(R'/J')+H_0$.
So we assume $c=\alpha$. By Lemma \ref{HFw2} we  have
$H(R/I)=H(R'/J')+H_\alpha$. Here  $J'=(\Ann G)\cap K[x,y,z]$ from
Lemma~\ref{alpha} has a
generator in degree $\alpha$, since by its definition \eqref{ealpha}
$\alpha $ is
the lowest degree for which $J'_i\nsubseteq \langle x,y\rangle
\cdot R'_{i-1}$. \par
 First, assume $\alpha < j/2$, when $H_\alpha=(0,1,\ldots
,1,2_\alpha,\ldots 2_{j-\alpha},1,\ldots
,1,0)$, from Definition \ref{halpha}.
We let
$H'=H(R/I)-H_0$, and we have
\begin{equation}\label{eH'}
H'=H(R'/J')+(0,\ldots ,0,1_\alpha ,1,\ldots
1_{j-\alpha},0,\ldots ).
\end{equation}
 Thus, to show \eqref{hfw2ii} here it would
suffice to show that
$H'$ of \eqref{eH'} is a height three Gorenstein sequence.
Assuming that the order of $J'$ is $\nu$, we have
\begin{align}\label{delta1}\Delta H'&=\Delta
H(R'/J')+(0,0,\ldots ,1_\alpha,0,\ldots
,0,-1_{j+1-\alpha},0,\ldots ),
\text{ and }\notag\\
\Delta H(R'/J') &=(1,2,3,\ldots ,\nu,t_\nu,\ldots ,-2,-1),
\end{align}
with $\nu\ge t_{\nu}\ge \ldots \ge t_{\lfloor j/2 \rfloor}$.
Furthermore, a result of A. Conca and G. Valla is that the maximum number
of degree-i
generators possible for any Gorenstein ideal
$J'$ of Hilbert function $H(R/J')$ is
\begin{equation*}
\max \{\nu_i\} =
\begin{cases} -(\Delta^2 H(R'/J')))_i=t_{i-1}-t_i, & \text{ when $ i\le j/2$
and  $ i\not=\nu$}\\ 
1-(\Delta^2 H(R'/J'))_\nu=1+\nu-t_\nu &\text { when 
$i=\nu$.}
\end{cases}
\end{equation*}
(see \cite{CV} or \cite[Theorem B.13]{IK}). Since $J'$ has a generator in
degree $\alpha$
it follows when $\alpha>\nu$ that $t_{\alpha -1}\ge t_\alpha +1$. Thus, for
$\alpha \ge \nu$ adding one in degree
$\alpha$ to the first difference
$(\Delta H(R'/J'))_{\le j/2}$ yields a sequence $\Delta H'$ as in
\eqref{delta1} that is still an
$O$-sequence: for height two this
condition is simply that the sequence $\Delta H'$ must rise to a maximum
value $\nu'$, then be nonincreasing. This implies that $H'$ is indeed
a height three Gorenstein sequence, and completes the proof
when
$\alpha
\le j/2$.\par
 Now assume that $c=\alpha$ and $\alpha > j/2$. Let
\begin{equation*}
H''=H(R'/J')+(0,\ldots ,0,-1_{j+1-\alpha} ,-1,\ldots ,
-1_{\alpha -1},0,\ldots ).
\end{equation*}
 Then we have in this case
$H(R/I)=H''+H_0$.  Thus, to show \eqref{hfw2ii} here it would
suffice to show that
$H''$ also is a height three Gorenstein sequence.   We have
\begin{equation}\label{delta2}\Delta H''=\Delta
H(R'/J')+(0,0,\ldots ,-1_{j+1-\alpha},0,\ldots
,0,1_{\alpha},0,\ldots ).
\end{equation}
That $J'$ has a generator in degree $\alpha> j/2$, implies that
$(\Delta^2(H(R'/J))_\alpha \le -1$, which is equivalent by the
symmetry of $\Delta ^2 (H(R'/J'))$ to $(\Delta ^2
(H(R'/J'))_{j+2-\alpha} \le -1$. This in turn implies
$(\Delta H(R'/J'))_{j+2-\alpha}< ( \Delta
H(R'/J'))_{j+1-\alpha}$. Thus, lowering $(\Delta
H(R'/J'))_{j+1-\alpha}
$ by 1 in degree $j+1-\alpha$ to obtain $\Delta
H''_{\le j/2}$ as in
\eqref{delta2} preserves the condition that
$(\Delta H'')_{\le j/2}$ is the Hilbert function of some height
two Artinian algebra. This completes the proof of the Theorem.
\end{proof}\par
  The following examples illustrate Lemma \ref{HFw2}. In particular we
explore how the Hilbert functions $H(R/I), H(R/J)$ change (recall
that 
$I=\Ann(F), J=\Ann(G)$)  as we alter the coefficient of $Z^{[j]}$
in $F,G$. Here there is a marked difference for the cases $\alpha (J)\le
j/2$, and $\alpha (J)> j/2$. The subsequent Lemma
\ref{lambda} explains some of the observations.
\begin{example}\label{anoninv}
Letting $G=X^{[4]}Z^{[2]}-X^{[4[}YZ, F=G+WZ^{[5]}$,
we have $J=\Ann(G)=(w,yz+z^2, y^2, x^5)$, so $\alpha (J)=2$, and
 $I=\Ann(F)=(
w^2,wx, wy, y^2 ,yz^2, xyz+xz^2, x^4y+wz^4, x^5, z^6)$. Also
$H(R/J)=(1,3,4,4,4,3,1)
$, and
\begin{equation*} H(R/I)=(1,4,6,6,6,4,1)=H(R/J)+H_2.
\end{equation*}
Changing $G$ by adding a $Z^{[6]}$ term, we  have
$G_1=X^{[4]}Z^{[2]}-X^{[4]}YZ+Z^{[6]},
F_1=G_1+WZ^{[5]}$, $J(1)=\Ann(G_1)=(w,y^2, yz^2,xyz+xz^2,x^4y+z^5,
x^5)$,
so $\alpha (J(1))=5$, and $I(1)=\Ann(F_1)=(w^2,wx, wy, y^2 ,yz^2,
xyz+xz^2, x^4y+wz^4, x^5, wz^5-z^6)$. Also
$H(R/J(1))=(1,3,5,5,5,3,1)
$, and
\begin{equation*} H(R/I(1))=(1,4,6,6,6,4,1)=H(R/J(1))+H_0.
\end{equation*}
\end{example}
\begin{example} In this example, we chose
$G=(Z+X)^{[6]}+(Z+2X)^{[6]}+(Z+Y)^{[6]}+(Z+2Y)^{[6]}+
(Z+X+Y)^{[6]}+(Z+2X+2Y)^{[6]}$,
the sum of 6 divided powers, and let $J=\Ann(G)$. Then $H(R/J)$
has the expected value $H(R/J)=(1,3,6,6,6,3,1)$ (see \cite{IK}), and $\alpha
(J)=3$. From
Lemma \ref{HFw2}, letting $
I=\Ann(F), F=G+WZ^{[5]}$ we have 
\begin{equation*}
H(R/I)=H(R/J)+H_3=(1,4,7,8,7,4,1).
\end{equation*}
Here
\begin{align*}
I=&(w^2,wx,wy, y^3-3y^2z+2yz^2, x^2y-xy^2 ,x^3-3x^2z+2xz^2 ,\\
&51xy^2z-18x^2z^2-99xyz^2-18y^2z^2-12wz^3+34xz^3+34yz^3, 5y^2z^3+4wz^4-9yz^4,
yz^5-z^6).
\end{align*}
\par
Omitting the pure $Z^{[6]}$ term from $G$ and $F$, to obtain $G_1,F_1$ we have
$H(R/\Ann  (G_1))=(1,3,6,7,6,3,1)$,  $\alpha(\Ann(G_1))=4$ and
\begin{equation*} 
H(R/\Ann F_1)=H(R/I)=H(R/\Ann (G_1))+H_0.
\end{equation*}
 This example
shows that it is not the inclusion of a $Z^{[6]}$ term in $G$ that keys
the simpler case
$H(R/I)=H(R/J)+H_0$.
The Hilbert function $H(R/I)$ is always invariant under a change
in the $Z^{[j]}$ term of
$F$: this follows from $z^i\circ F=
WZ^{[j-1-i]}+z^i\circ G$, linearly disjoint from $\langle R_i\mod z^i\rangle
\circ F$.
\end{example}
\begin{example} When $j=8$,
$G=X^{[3]}Y^{[5]}+X^{[2]}Y^{[4]}Z^{[2]}+Y^{[5]}Z^{[3]}$, then
\begin{equation*}
J=\Ann G=
(w,x^3-z^3, z^4 ,xz^3, x^2z^2-yz^3, y^6, xy^5+x^2y^3z-y^4z^2,
x^2y^4-y^5z-xy^3z^2).
\end{equation*}
We have $\alpha (G)=3$, $H(R/\Ann (G))=(1,3,6,9,9,9,6,3,1)$, and $I
=\Ann(F), F=G+WZ^{[7]}$,
satisfies 
\begin{equation*}
H(R/I)=(1,4,7,11,11,11,7,4,1)=H(R/\Ann (G))+H_3.
\end{equation*}
Here
\begin{align*}I=(w^2,wx,wy,&
 xz^3, x^2z^2-yz^3, x^3z, x^3y-yz^3, x^4, y^5z-wz^5, y^6,\\&
xy^5+x^2y^3z-y^4z^2,
x^2y^4-xy^3z^2-wz^5  ,z^8).
\end{align*}
Adding a $Z^{[8]}$ term to $G$ to form $G_1$ leads to $ J(1)=\Ann (G_1)$ with
$\ \alpha (J(1)) =6$ and $F_1, I(1)=\Ann (F_1)$ satisfying
\begin{equation*}
H(R/I(1))=H(R/I)=H(R/J(1))+H_0.
\end{equation*}
\end{example}
It might be thought from the previous examples, that adding $\lambda
Z^{[j]}$ with
$\lambda$ generically chosen, will ``improve'' $G$ to a $G_\lambda$ such that
$J(\lambda)=\Ann G_\lambda$ and $I_\lambda=\Ann F_\lambda, F_\lambda
=G_\lambda+WZ^{[j-1]}$ will satisfy
$H(R/I_\lambda)=H(R/J_\lambda) +H_0$. This change would indeed be an
improvement, since when $H(R/I)=H(R/J)+H_0$ the minimal resolutions of
the ideals $I,J$ appear to be closer than they are when
$H(R/I)=H(R/J)+H_\alpha$.
In the next Lemma we show that this ``improvement'' must occur when $\alpha
(J)\le j/2$, but can occur either never, or for a single value of
$\lambda$ when
$\alpha (J) > j/2$. We suppose that
$\lambda
\in K$.
\begin{lemma}\label{lambda} Let $J=\Ann(G), I=\Ann(F), F=G+WZ^{[j-1]}$ be
such that
$I$ satisfies
\eqref{ebase}, and define $G_\lambda=G+\lambda Z^{[j]}$, $F_\lambda = F+\lambda
Z^{[j]}$,
$J(\lambda)=\Ann (G_\lambda), I(\lambda)=\Ann(F_\lambda)$. Then we have
\begin{enumerate}[i.]
\item\label{lambdai} $\left( I\cap J\right) +m^j=\left(
I(\lambda)
\cap J(\lambda)\right) +m^j$ and $(I\cap J)_j$ differs
from $(I(\lambda )
\cap J(\lambda ))_j$ by replacing $z^j-u,u\in J\cap ((x,y)\cap
K[x,y,z])$ by $z^j-u', u'\in J(\lambda )\cap ((x,y)\cap
K[x,y,z])$.
\item\label{lambdaii} $H(R/I)=H(R/I(\lambda)), $ and $H(R/(I\cap
J))=H\left( R/(I(\lambda)\cap J(\lambda))\right)$;
\item\label{lambdaiii} If $\alpha (J)\le j/2$ and $\lambda \ne 0$
then $\alpha (J(\lambda))=j+1-\alpha (J)$, and \begin{equation*}
H(R/J(\lambda ))=
H(R/J)+(0,\ldots,0_{\alpha -1},1_\alpha ,1,\ldots
,1_{j-\alpha},0_{j+1-\alpha } ,\ldots ,0_j).
\end{equation*}
 In this case
$H(R/I(\lambda ))=H(R/J(\lambda ))+H_0$.
\item\label{lambdaiv}  Let $\alpha (J)>j/2$
then $\ \alpha (J(\lambda))=\alpha(J)$ or $\alpha (J(\lambda))=j+1-\alpha
(J)$. In the former case
$H(R/J(\lambda ))=H(R/J)$. The latter case may occur for at most a single
value $\lambda _0$; if it occurs, then for $\lambda =\lambda
_0,\alpha=\alpha(J)$,
\begin{equation*}
H(R/J(\lambda _0))=H(R/J)-(0,\ldots,0_{j-\alpha },1_{j+1-\alpha}
,1,\ldots ,1_{\alpha -1},0_{\alpha } ,\ldots ,0_j).
\end{equation*}
\begin{enumerate}[a.]
\item\label{lambdaiva} If
$H(R/I)=H(R/J)+H_\alpha$ then $\alpha
(J(\lambda ))=\alpha (J)$ and $H(R/J(\lambda))=H(R/J)$.
\item\label{lambdaivb} If
$H(R/I)=H(R/J)+H_0$, then for all values of $\lambda$ except
possibly a single value $\lambda _0\not= 0$ we have
$\alpha(J(\lambda))=\alpha (J)$ and $H(R/J(\lambda ))=H(R/J).$
\end{enumerate}
\end{enumerate}
\begin{proof}Since for $i\le j-1, Z^{[i]}=wz^{j-1-i}\circ
WZ^{[j-1]}$, $(I\cap J)_i=(I(\lambda ) \cap J(\lambda ))_i$ for
$i\le j-1$. The second statement in \eqref{lambdai} is evident.
The first claim in \eqref{lambdaii} follows since the two ideals
$I, I(\lambda)$ are isomorphic, under a change of variables. The
second claim in \eqref{lambdaii} follows from \eqref{lambdai}. \par
Suppose that $\alpha \le j/2$ and $\lambda\ne 0$, and that $h=z^\alpha-g,g\in
(x,y)\cdot K[x,y,z] \in J$. Then for $0\le u \le j-2\alpha $ we have
\begin{equation*}
(z^u
h) \circ (G+\lambda Z^{[j]})= z^uh\circ (\lambda Z^{[j]})=\lambda
Z^{[j-\alpha-u]}.
\end{equation*}
It follows that for $i\le  j-\alpha$, $Z^{[i]}\in R\circ G(\lambda)$. This
implies
that for $\alpha\le i \le j-\alpha$, we have $ H(R/J(\lambda )_i=H(R/J)_i+1$,
since by Lemma \ref{HFw2} $Z^{[i]}\ne R\circ G $ for $i\ge \alpha (J)$. The
claims
in
\eqref{lambdaiii} now follow from the
symmetry of
$H(R/J(\lambda)), H(R/J)$ and hence of $H(R/J(\lambda))-H(R/J)$.\par
Suppose that $\alpha (J)> j/2$. The symmetry of $H(R/J(\lambda))-H(R/J)$
and Lemma \ref{HFw1} \eqref{HFw1ii} show the first claim concerning $\alpha
(J(\lambda))$ in
\eqref{lambdaiv}. This and \eqref{lambdaii} show \eqref{lambdaiva}.
The same symmetry, and \eqref{lambdaiii} also prove \eqref{lambdaivb}, and
completes the proof of \eqref{lambdaiv} that the exceptional case may occur for
at most a single value $\lambda _0$.
\end{proof}
\end{lemma}
\begin{example} Letting 
$G=X^{[3]}Z^{[3]}-Y^{[4]}X^{[2]}+Y^{[2]}Z^{[4]}+XY^{[2]}Z^{[3]}+Z^{[6]},
F=G+WZ^{[5]}$ we have $J=\Ann(G)=(x^3+x^2z-y^2z, y^3z,
y^4-x^2z^2+y^2z^2+xz^3-z^4, xy^2z+x^2z^2-y^2z^2, xy^3-xyz^2+yz^3, x^2yz,
x^2y^2+x^2z^2-y^2z^2+z^4)$, $\alpha (J)=4$, and $H(R/J)=(1,3,6,9,6,3,1)$.
Then
\begin{align*} I&=(wy, wx, w^2, x^2z-y^2z+wz^2, x^3-wz^2 ,y^3z,
xy^3-xyz^2+yz^3, x^2y^2+y^4+xz^3, wz^5-z^6), \text{ and }\\
&\qquad \qquad H(R/I)=(1,4,7,9,7,4,1)=H(R/J)+(0,1,1,0,1,1,0)=H(R/J)+H_4.
\end{align*}
This is an example of Lemma \ref{lambda}\eqref{lambdaiva} where
$H(R/J(\lambda))=H(R/J)$ for every $\lambda$.
\end{example}
\medskip

\section{Hilbert functions $H=(1,4,7,h,\ldots ,4,1)$}\label{jarb}
We now consider Gorenstein sequences --- Hilbert functions of Artinian
Gorenstein algebras, so
symmetric about $j/2$ --- having the form
\begin{equation}\label{genhfe}
H=(1,4,7,h,b,\ldots ,4,1),
\end{equation}
 of
any socle degree
$j\ge 6$ for any possible $b$.  We show in Theorem \ref{nonempty2} that each such
Gorenstein sequence must satisfy the \emph{SI
condition} that
$\Delta H_{\le j/2}$ is an O-sequence. This
condition was shown by R. Stanley, and
by D.~Buchsbaum and D. Eisenbud to characterize Gorenstein sequences of height
three (see \cite{BE,St,Hari2}).
When a Gorenstein sequence $H$ satisfies this condition
 we can construct Artinian Gorenstein algebras, elements of $\PGOR(H)$, as
quotients of the
coordinate ring of suitable punctual schemes, and we have good control
over their Betti numbers (Lemma \ref{sevcomp2AB}, Corollary \ref{relGor}). In
particular, when
$H=(1,4,7,h,\ldots )$ satisfies the {SI} condition and $7\le h\le 10$ we
may choose $A\in
\PGOR (H)$ such that
$I_2$ has only two linear relations: thus
$A\notin \overline{{\mathfrak{C}}(H)}$, the locus where $I_2\cong
\langle wx,wy,wz\rangle$, implying for most such Hilbert functions $H$ that
$\PGOR(H)$ has at least two irreducible components (Theorems
\ref{sevcompt}, \ref{sevcomp2}). \par Our first result is relevant also to the
open question of whether all height four Gorenstein sequences satisfy the
SI condition. Despite our positive result we doubt that this is true in
general (see Remark
\ref{SIconj}). We now set some notation. When $H$ is clear we usually write $h_i$
for
$H_i$ below.. We set $\Delta
H_i=h_i-h_{i-1}$. By
$H_{i,i+1}$ we mean $(h_i,h_{i+1})$.
 Given a Hilbert function $H_\Z$, we define $\Sym(H_\Z,j)$ as the
symmetrization of
$(H_\Z)_{\le j/2}$ about $j/2$:
\begin{equation} \Sym(H_\Z,j)_i=
\begin{cases}
(H_\Z)_i \text{ if $i\le j/2$}\\
(H_\Z)_{j-i} \text{ if $i> j/2$}
\end{cases}.
\end{equation}
\begin{lemma}\label{nonempty2a}
 Let $j\ge 6$ and suppose that the Gorenstein sequence $H$ of socle degree
$j$ satisfies
\eqref{genhfe}. Then
$7\le h\le 11$. If $j\ge 7$, then the minimum value of $b=H_4$ that can
occur is
$b=h$, and the maximum values of
$b$ that can occur in \eqref{genhfe} are
\begin{equation}\label{nonemptyeq}
\begin{array}{r|ccccc}
 h & 7 & 8 & 9 & 10 & 11\\
b_{\max}  & 7&9&11&13&16
\end{array}
\end{equation}
Equivalently, a Gorenstein sequence $H$ satisfying \eqref{genhfe} must
satisfy $\Delta H_{\le 4}$ is an
$O$-sequence.
Also, each initial sequence $(1,4,7,h,b)$ satisfying $7\le h\le 11$ and
$h\le b\le b_\max $ occurs for
$j=8$. \par
Finally, if $H$ satifies \eqref{genhfe} and $j\ge 6,h\le 10$ then $\Delta
H_{\le j/2} $ is an $O$-sequence if and only if its subsequence $\Delta
H_{1\le i\le
j/2}=(3,3,h-7,b-h,\ldots )$ is both nonnegative and nonincreasing.
\end{lemma}
\begin{proof}
We showed $7\le h\le 11$ in Corollary \ref{nonempty1}.  We
now show the upper bounds
$b\le b_\max$ of \eqref{nonemptyeq}. When
$h=11$, the upper bound of \eqref{nonemptyeq} is just the Macaulay upper
bound. When $h=10$, the
impossibility of
$(h,b)=(10,15)$ follows from Corollary \ref{oseq}. The impossibility of
$(h,b)=(10,14)$
follows from two considerations.
First, by Theorem~\ref{V,W}
\eqref{V,Wiiia} and Theorem \ref{hfw2}\eqref{hfw2ii} $I_2$ cannot be
$\PGL(3)$-isomorphic
to $\langle wx,wy,wz\rangle
$ or
$\langle w^2,wx,wz\rangle$, as $H'=H-(0,1,1,\ldots ,1,0)=(1,3,6,9,13,\ldots
)$ is not a
height three Gorenstein sequence, since $\Delta H'_{\le
j/2}=(1,2,3,3,4,\ldots )$ is not an O-sequence in
two variables  \cite{BE,D}.  Thus $I_2$ cannot have a common factor, so has
two linear relations. By
Lemma
\ref{netsq}
\eqref{netsqi} $I_2$ has a basis given by the $2\times 2$ minors of a
$2\times 3$ matrix; since $I_2$
has no common factor, the quotient $R/(I_2)$ has height two, $I_2$ is
determinantal
and has the usual determinantal minimal resolution. In
particular we have
$H(R/(I_2))_i=3i+1$, for all $i\ge 0$, so as before
$H(R/I)_4\le H(R/(I_2))_4=13$. \par
When $h=8$ or $9$ the upper
bound of
\eqref{nonemptyeq} is one less than the Macaulay upper bound. The
impossibility of the Macaulay upper bound
for $H(R/I)_{3,4}$ in the cases $h=8,9$ follow from
Lemma~\ref{Hrest}\eqref{Hrestiii}. When $h=7$, 
the upper bound $b\le 7$ is shown in the $h=7$ case of the proof of
Theorem \ref{nonempty2} below. This completes the proof of the upper
bounds $b\le b_\max$ of \eqref{nonemptyeq}.
\par We next show the lower bound on $b$: when $j\ge 7$, then $b\ge h$.
Evidently, when $j=7$, the symmetry
of
$H$ implies
$b=h$, so we may assume
$j\ge 8$. The symmetry of $H$ implies $(H_{j-4,j-3})=(b,h)$. The   Macaulay
Theorem
\ref{MacGo}
\eqref{Macgrowth} applied to  $(H_{j-4,j-3})$ eliminates all
triples
$(j,h,b)$ where
$b\le h-2$ except the triple $(j,h,b)=(8,5,4)$. For this triple
$H_{4,5}=(b,h)=(9,11)$ is extremal
growth as
$9^{(4)}=11$; then we have a contradiction by Corollary \ref{oseq}.\par
We now assume $j\ge 8$ and $b=h-1$. We have $h\neq 11$ by Theorem \ref{V,W}
\eqref{V,Wiiia} and Theorem \ref{hfw2}. Since in Macaulay's inequality of
Theorem
\ref{MacGo}\eqref{Macgrowth} $b^{(d)}=b$ when $b\le d$, and $h_{j-4}=b,
h_{j-3}=b+1$ we must have
$b>j-4$, so $h\ge j-2$. Except for the triples $(j,h,b)=(8,10,9)$ or
$(8,11,10)$, then
$H_{j-4,j-3}$ has extremal Macaulay growth, a contradiction by Corollary
\ref{oseq}. The second
triple has
$h=11$, already ruled out. The first triple occurs only for
$H=(1,4,7,10,9,10,7,4,1)$
where $\Delta^4 H_6=-12$;  by symmetry of the minimal resolution of $R/I$,
the number of degree six
generators of $I$ satisfies
$\nu_6(I)\ge 6$, implying that
$H(R/(I_5))_{5,6}=(10,13)$, contradicting the Macaulay bound which requires
$H(R/(I_5))_6\le
10^{(6)}=11$. This completes the proof of the lower bound on $b$, that
$h\le b$ in \eqref{genhfe}.\par
It is easy to see that these bounds are just the condition that $\Delta
H_{\le 4}$ be an $O$-sequence,
as claimed.
\par That
each extremal pair
$(h,b)$ satisfying $h\le b\le b_\max$ from \eqref{nonemptyeq} occurs in
socle degree 8 can be shown by
choosing the ring $A$ to be a general enough socle-degree 8 Artinian
Gorenstein quotient of the
coordinate ring of any smooth punctual scheme of degree
$b$, having Hilbert function $H_\Z=(1,4,7,h,b,b,\ldots )$. Since $b\ge h$,
$\Delta H_\Z$ is an
$O$-sequence and there are Artinian algebras of Hilbert function $\Delta
H_\Z$; then there is a smooth punctual
schemes of Hilbert function $H_\Z$, by the result of P. Maroscia
\cite{Mar,GMR,MiN}). That the general socle-degree $j$ GA quotient of
$\Gamma (\Z,{\mathcal{O}_\Z})$ has
the expected symmetrized Hilbert function $H=\Sym(H_\Z,j)$ satisfying
$(\Sym(H_\Z,j))_i=(H_\Z)_i$ for
$i\le j/2$, is well known: see
\cite{Bj1,MiN}\cite[Lemma 6.1]{IK}.\par
The last statement of Lemma \ref{nonempty2a} that $j\ge 6,h\le 10$ and
$\Delta H_{\le j/2}$ an
$O$-sequence is equivalent to
$\Delta H_{2\le i\le j/2}$ being non-negative and non-increasing, follows
from $\Delta H=(1,3,3,h-7,\ldots )$, with $h-7\le 3$: by Macaulay's
inequality Theorem \ref{MacGo}\eqref{Macgrowth}, we have for any
$O$-sequence $T$ that $t_i\le i$ implies $t_{i+1}\le i$.
\end{proof}
\begin{theorem}\label{nonempty2}
Every Gorenstein sequence $H$ beginning $H=(1,4,7,\ldots )$
satisfies the condition,
$\Delta H_{\le j/2}$ is an O-sequence. \par
\end{theorem}
\begin{proof} We assume $H=H(R/I)$ for an Artinian Gorenstein quotient
$R/I$ satisfies
\eqref{genhfe} that $H=(1,4,7,h,b,\ldots )$ and consider each value of
$h$ in turn.  We show that each occurring sequence $H$ satisfies the
criterion from Lemma
\ref{nonempty2a} for
$\Delta H$ to be an
$O$-sequence.\par {\bf Case}
$h=7$. We have $H(R/I)_{j-3,j-2}=
H(R/I)_{2,3}=(7,7)$; if $j\ge 10$ then $H$ is extremal in
degrees j-3 to j-2, and we have that
$\Z=\Proj(I_{j-3})$ is a degree-7 punctual scheme satisfying by Lemma
\ref{upbd}
$H(_\Z)_i=7$ for all $i\ge 3$: by Corollary \ref{include}, we have
$H(R/I)_i=7$ for $3\le i\le
j-2$. So we may assume that $j=8$ or $9$.  We have
$b\le 7^{(3)}= 9$. Should
$b=9$ then $\Proj(R/(I_3))$ would
define a degree-2 curve of genus zero and regularity two, so its  Hilbert
function would satisfy
$H(R/(I_\Z))_2\le 5$, by Corollary \ref{include} contradicting
$H(R/I)_2=7$. We now suppose that $h=7,b=8$, and suppose the socle degree
 $j=8$ or $9$. When $j=8$, $H=(1,4,7,7,8,7,7,4,1)$, since $\Delta^4
H_5=-7$, the ideal
$I$ has
 $\nu_5$ generators (first syzygies)  and $\mu_5$ third syzygies in degree
5, with $7\le\nu_5+\mu_5$; by symmetry of the minimal resolution
$\nu_7=\mu_5$ and
$\mu_7=\nu_5$; Thus we have either $\nu_5\ge 3$ or $\nu_7\ge
4$; but $\nu_5\le 2$ and $\nu_7\le
4$ by Macaulay's Theorem \ref{MacGo}. If $\nu_7=4$ then the ideal $(I_{\le
6})$ would  satisfy $H(R/(I_
6))_{6,7}=7,8$ of extremal growth, a contradiction with $\Delta H_3=0$, by
Corollary \ref{oseq} and
Lemma
\ref{nonempty2a}. For $j=9$ we would have similarly $\Delta^4 H_5=-6$, so
$\nu_5+\nu_8\ge 6$, but
$\nu_5\le 2$, and when $\nu_8=4$ we'd have $H(R/(I_
6))_{7,8}=7,8$, and a similar contradiction.
We have shown that a Gorenstein sequence beginning
$(1,4,7,7)$ continues with a subsequence of $7's$ followed by $(4,1)$.
\par {\bf Case} $h=8$. Macaulay extremality shows $h_i\le i+5$ and $\Delta
H_{i+1}\le 1$ for $i\ge
3$. Suppose by way of contradiction that
$\Delta H_i<0$, for some  $i\le j/2$ (this is equivalent to $H$ being
non-unimodal).  Letting
$i'=j-i$, we have by the symmetry of
$H$ that
$h_{i'+1}=h_{i'}+1=h_{i-1}\le i-1+5\le i'+4$; it follows from Theorem
\ref{MacGo} that either this
is impossible (when $h_{i'}\le i'$) or
$H$ is extremal in degrees $i'$ to $i'+1$, a contradiction by Corollary
\ref{oseq} and Lemma
\ref{nonempty2a}. Now suppose
$4\le i<k,\Delta H_{i}=0$ but
$\Delta H_{k}=1$. Then $H_k=(H_{k-1})^{(k-1)}$, and we have a contradiction
by Corollary
\ref{oseq} and Lemma \ref{nonempty2a}. It follows that
$H$ satisfies,
$\Delta H_i, 2\le i\le j/2$ is nonnegative
and nonincreasing, thus $\Delta H_{\le j/2}$ is an $O$-sequence.\par
{\bf Case $h=9$}. Lemma \ref{Hrest}\eqref{Hrestiii} implies that $h_4\le
11$; applying  Macaulay extremality
inductively we have for $i\ge 4$ that
$h_i\le 2i+3$ and $\Delta H_i\le 2$. Suppose by way of contradiction that
$\Delta H_i<0$, for some  $i\le j/2$; then $h_i=i+a$ with $a\le i$. We now
 use
the symmetry of $H$ about $j/2$. Letting
$i'=j-i$, we have
$h_{i'}=i+a=i'+a',a'=a-(i'-i)$; since $a'<a<i'$ we must have $h_{i'}\le
2i'$ whence $h_{i'+1}\le
h_{i'}+1$ by the Macaulay Theorem
\ref{MacGo}\eqref{Macgrowth}, so $\Delta H_{i'+1}=-\Delta H_i=1$, and
$h_{i'+1}=h_{i'}+1$ is extremal, a contradiction by Corollary \ref{oseq}
and Lemma
\ref{nonempty2a}.\par Now suppose that for some
$i\le j/2$ we have $\Delta H_{i-1}=1$, but
$\Delta H_i=2$: then by Theorem~\ref{MacGo} we would have $\Proj(R/(I_i))$
defines a degree 2
curve union some points,  of Hilbert polynomial $2t+a,a\le 2$, of
regularity degree at most 3 by  Corollary \ref{regdeg}, hence by Lemma
\ref{upbd} and
Corollary~\ref{include} we would have $h_3\le 8$, a contradiction. \par
Finally, suppose that for
some $i\le j/2$ we have
$\Delta H_{i-1}=0$ but
$\Delta H_i>0$. By Corollary \ref{oseq} we have $\Delta H_i\not= 2$, so
$\Delta H_i=1$. If
also there is a previous $u, 4\le u\le i-2$ with $\Delta H_u<2$ then
$h_i\le 2i$, implying
that $H_i=(H_{i-1})^{(i-1)}$, a contradiction by Corollary \ref{oseq}. Thus
to complete
the case $h=9$, we need only
consider sequences
\begin{equation}\label{h=9lastcase}
H=(1,4,7,9,\ldots ,h_u=2u+3,\ldots ,h_{i-2}=h_{i-1}=2i-1, h_i=2i,\ldots ,7,4,1)
\end{equation}
 with possible consecutive
repetition of the maximum value $2i$. We have $\Delta^4
H_{i+1}=-5$ if $h_{i+1}=h_i$, and
$-6$ if $j=2i$ so $h_{i+1}=h_i-1$. In either case, we obtain
$\nu_{i+1}+\nu_{j+3-i}\ge 5$. This is
impossible since on the one hand $\nu_{j+3-i}\ge 3$ would
imply that
$H(R/(I_{j+2-i}))_i=h_{i-2}=2i-1, H(R/(I_i))_{j+3-i}=h_{i-3}+3=2i-3+3=2i$,
which is extremal growth of $H$, a contradiction by Corollary \ref{oseq}.
On the other hand if
$\nu_{i+1}\ge 1$ when
$h_{i+1}=h_i$, or if
$\nu_{i+1}\ge 2$ when
$h_{i+1}=h_i-1$ we would have
$H(R/I)_i=2i,H(R/(I_i))_{i+1}=2i+1 $ implying
extremal growth, a contradiction with \eqref{h=9lastcase} by Corollary
\ref{oseq}.
This completes the proof that $\Delta H$ is an
$O$-sequence when $h=9$.\par {\bf Case }$h=10$. By Lemma \ref{nonempty2a}
$h_4\le 13$; also
when $I_2$ has a common factor Theorems \ref{V,W}\eqref{V,Wiiia} and
Theorem
\ref{hfw2} show that $\Delta H_{\le j/2}$ is an $O$-sequence.
We suppose henceforth in our analysis of $h=10$ that
$I_2$ does not have a common factor. Then by Lemma
\ref{netsq}\eqref{netsqi}
$I_2$ defines a rational normal curve, satisfying
$H(R/(I_2))_t=3t+1$ for all $t\ge 0$.
 Notice also that if $H(R/I)_t\le
3t-1$, and $t\ge 4$, then the Macaulay inequality Theorem
\ref{MacGo}\eqref{Macgrowth} implies $\Delta H(R/I)_{i+1}\le
2$. We next rule out various perturbations in the Hilbert
function sequence.
\par First, $\Delta H_{i+1}\le -2$ for some $i<j/2$ is impossible from the
Macaulay bound and the symmetry of
$H$. We would have $\Delta H_{i'+1}\ge 2$ for
$i'=j-i-1\ge i+1$;
then letting
$h_i=3i+1-e, e\ge 0$ we have
$h_{i'}=h_{i+1}\le h_i-2= 3i-(e+1)=2i+(i-e-1)=2i'+ b, b\le i-e-3$; thus, the
Macaulay bound here implies $\Delta H_{i'+1}\le 2$, so there is equality
$\Delta H_{i'+1}=2$,
a contradiction by Corollary \ref{oseq}. Also $\Delta H_{i+1}=-1 $ for some
$i<j/2$,
and
$j> 5i+e$, is impossible by a similar calculation that $\Delta H_{i'+1}=1$
the maximum possible, again a contradiction by Corollary
\ref{oseq}. \par
 Suppose $\Delta H_{i+1}=-1 $ with $i\le j/2-1$ and no restriction on $j$;
suppose that $i$ is the
maximum such integer. Letting
$ c=h_{i+1}$ we write the consecutive subsequence
$(h_{i-1},\ldots , h_{i+3})$ as
\begin{equation}\label{h=10,dec1}
(a+c,1+c,c,1-\alpha+c, b+c).
\end{equation}
Then $\nu_{i+3}(I)+\nu_{j+5-i}(I)\ge -\Delta^4 H_{i+3}=-\Delta^4
H_{j+5-i}=8-a-b-4\alpha$. We have
\begin{align*}H(R/(I_{i+2}))_{i+2,i+3}&=(1-\alpha+c,b+c+\nu_{i+3})\hbox{
and }\\
H(R/(I_{j+4-i}))_{j+4-i,j+5-i}&=(1+c,a+c+\nu_{j+5-i}).
\end{align*}
Thus the sum $\delta +\delta', \delta =\Delta H(R/(I_{i+2}))_{i+3},\delta
'=H(R/(I_{j+4-i}))_{j+5-i}$ satisfies
\begin{equation*}
\delta +\delta '=(b+\nu_{i+3}+\alpha-1)+(a+\nu_{j+5-i}-1)\ge 6-3\alpha .
\end{equation*}
 So
if $\alpha\le 1$ at least one of $\delta ,\delta '$ is two, and the
corresponding
Hilbert function has extremal growth of two, a contradiction by Corollary
\ref{oseq}. If $\alpha
=2$, then $i+1\le j/2-1$ (by the symmetry of $H$), and $\Delta H_{i+2}=-1$,
contradicting the
assumption on
$i$; and
$\alpha\ge 3$ has already been ruled out. We have shown $\Delta H_{i+1}=-1$
for $i\le j/2-1$
is impossible.\par We cannot have both $\Delta H_u\le 2$ and $\Delta
H_{i+1}=3$ for a pair $u,i$
satisfying
$u<i<j/2$, since then
$h_i\le 3i$. This is possible only if $h_i=3i$ and $h_{i+1}=h_i^{(i)}$,
 a contradiction by
Corollary~\ref{oseq}. We cannot have both $\Delta H_u\le 1$ and $\Delta
H_{i+1}=2$ for $u<i<j/2$,
 since then
$h_i= 3i-1-e,e\ge 0$, and $H_{i,i+1}$ is extremal, again a contradiction by
Corollary \ref{oseq}.
\par
Suppose that for some $i, 2\le i\le j/2-1$, we have $\Delta H_i=0$, but
$\Delta H_{i+1}=1$.
 Then, letting $ c=h_i$ the consecutive subsequence
$(h_{i-2},\ldots , h_{i+2})$ is
\begin{equation}\label{h=10,delta0,1}
(a+c,c,c,1+c, b+c).
\end{equation}
Then $\nu_{i+2}(I)+\nu_{j+6-i}(I)\ge -\Delta^4 h_{i+2}=-\Delta^4
h_{j+6-i}=4-(b+a)$. It follows that the sum
$\Delta H(R/(I_{i+2}))_{i+3}+\Delta
H(R/(I_{j+5-i}))_{j+6-i}=a+b-1+4-(a+b)=3$, hence one of the two
differences is at least two, which is here extremal growth, since
$H_{i+2}\le 3(i+2)$
and similarly $H_{j+5-i}\le 3(j+5-i)$. Then Corollary
\ref{oseq} implies a contradiction with
\eqref{h=10,delta0,1}.\par
This completes the proof in the case $h=10$.\par
 {\bf Case }$h=11$. In this case $I_2$ must have a common linear factor.
 Theorem \ref{V,W}
\eqref{V,Wiiia} for $I_2\cong \langle wx,wy,wz\rangle$ and Theorem
\ref{hfw2} for $I_2\cong \langle w^2,wy,wz\rangle$ show that
$H=H'+(0,1,1,\ldots ,1,0)$, which
implies that $\Delta H_{\le j/2}$ is an $O$-sequence. \par
This completes the proof of the Theorem.
\end{proof}\par
For $H$ satisfying \eqref{genhfe}, recall that we denote by $\mathfrak
{C}(H)\subset
\PGOR(H)$ the subfamily parametrizing ideals $I$ such that $I_2\cong
\mathfrak{V}=\langle
wx,wy,wz\rangle$, up to a coordinate change. By Theorem
\ref{V,W}\eqref{V,Wii} we have that
$\mathfrak {C}(H)$ is nonempty if and only if $\PGOR(H')$ is nonempty, where
$H'=(1,3,6,h-1,i-1,\ldots , 3,1)$.
\begin{cor}\label{nonempty3} Let $H=(1,4,7,\ldots )$. The following are
equivalent.
\begin{enumerate}[i.]
\item\label{nonempty3i} The sequence $H$ is a Gorenstein
sequence.
\item\label{nonempty3ii} The sequence $\Delta H_{\le j/2}$ is an O-sequence.
\item\label{nonempty3iii} The sequence $H'=H-(0,1,1,\ldots ,1,0)$ is a
height three Gorenstein sequence.
\item\label{nonempty3iv}
$\Delta H'_{\le j/2}$ is an O-sequence.
\item\label{nonempty3v}  $\Delta H'_{\le j/2}=(1,2,3,\ldots
,i+1,h_v,h_{v+1},\ldots )$ with $i+1\ge h_v\ge
h_{v+1}\ge \ldots $.
\end{enumerate}
Under this assumption, the subfamily
$\mathfrak {C}(H)\subset
\PGOR(H)$ is always nonempty.
\end{cor}
\begin{proof}
That \eqref{nonempty3i} is equivalent to \eqref{nonempty3ii} is Theorem
\ref{nonempty2}. That
\eqref{nonempty3ii} is equivalent to
\eqref{nonempty3iv} is immediate from the last statement of Lemma
\ref{nonempty2a}, and an easy verification when $H=(1,4,7,11,\ldots )$. That
\eqref{nonempty3iii} is equivalent
to \eqref{nonempty3iv} follows from the Buchsbaum-Eisenbud structure theorem
\cite{BE,St}. That specific criterion
\eqref{nonempty3iv} is equivalent to \eqref{nonempty3v} is well known ---
see for example
\cite[Theorem 5.25,Corollary C6]{IKl}. That
$\mathfrak {C}(H)$ is always nonempty when $H$ satisfies these
conditions follows from Theorem
\ref{V,W} and
\eqref{nonempty3iii}.
\end{proof}\par
The following result handles height four Gorenstein sequences below those
considered in Theorem
\ref{nonempty2}.
\begin{proposition}\label{aless7} A symmetric sequence $H=(1,4,a,\ldots
,4,1), a\le 6$ of socle degree
$j$ is a Gorenstein sequence if and only if
$\Delta H_{\le j/2}$ is an O-sequence, or, equivalently, if
$\Delta H_{\le j/2}$ is nonincreasing once it does not increase. The values
$a=2,3$ cannot occur.
\end{proposition}
\begin{proof} When $a=6$, then $h_3=10$, the maximum under Macaulay's
theorem, would imply $h_1=3$, by Corollary \ref{oseq}. Assume $H=H(A)$ for
an Artinian
Gorenstein $A=R/I$ and let $\alpha_i$ denote the
number of relations (first syzygies) in degree $i$. When
$a=6$ and
$h_3=9,h_4=b$, then the fourth differences of $H$ satisfy
$\Delta^4(H)_4=\Delta^4(H)_j=b-15$, so by the symmetry of the
minimal resolution of
$A$ we have $\alpha_4+\alpha_j\ge 15-b$.  Since $\Delta
H(R/(I_{j-1})_j=\alpha-3$ and $j-1\ge 5$,
the Macaulay bound implies that growth from $h_{j-1}=4$
to $H(R/(I_{j-1}))_j=1+\alpha$ would be maximal when $\alpha_j=3$. But
$\alpha_j=3$, is impossible by Corollary \ref{oseq}.
However,  $\alpha_j\le 2$, implies $\alpha_4\ge
13-b$; thus $H(R/(I_3))_4\ge
b+13-b=13$, contradicting the Macaulay bound of
$9^{(3)}=12$. We have shown
$H=(1,4,6,9,\ldots )$ to be impossible. Establishing the result for
$H=(1,4,6,b,\ldots )$ with $b\le 8$ is
relatively simple, requiring only Theorem \ref{MacGo} and Corollary
\ref{include} without using the
symmetry of the minimal resolution: we leave this to the reader. \par
When $a=5$, then the Macaulay bound gives $h_3\le 7$; and
$H=(1,4,5,7,b,\ldots )$ is not possible by Corollary \ref{oseq}.
\par
The remaining cases are simpler, and we leave them as an exercise. Note
that $a=2,3$ are impossible,
since by the symmetry of $H$, we would have $h_{j-2}=a$ and $h_{j-1}=4$:
however, the
Macaulay bound gives $a^{(j-2)}\le a$ when $a\le j-2$, and here $j-2\ge 4$.
\end{proof}
\begin{remark}\label{SIconj} {\sc Do height four Gorenstein sequences
satisfy $\Delta H_{\le j/2}$ is
an $O$-sequence?} The height four Gorenstein sequences of the form
$H=(1,4,7,\ldots )$ are probably close to an upper
bound of those which may be shown to satisfy the condition $\Delta
H_{\le j/2}$ is an
$O$-sequence, by the kind of arguments we have used for Theorem
\ref{nonempty2}. Notice that we were not able to  rule out the nonoccuring,
sequence
$H=(1,4,7,10,14,10,7,4,1)$ by  a simple application of Macaulay bounds
and the Gotzmann
method of Lemma \ref{MacGo}, together with calculation of $\Delta^4 H$.
Rather, we needed to use
Lemma \ref{netsq}, which involves the twisted cubic. Likewise, in
proving other parts of Theorem~\ref{nonempty2}, we use at times detailed
information about low degree curves in
$\mathbb P^3$.
\par Thus we are inclined to conjecture that there are height four
Gorenstein sequences that do not
satisfy the condition that $\Delta H_{\le j/2}$ is an
$O$-sequence.
\end{remark}
 Recall that we denote by
$\nu_i(J)$ the number of degree-$i$ generators of the ideal $J$. The next
result
follows from Theorems \ref{sevcomp} and \ref{nonempty2}. Recall that the
socle degree of $H$ is the highest $j$
such that
$h_j\not=0$.
\begin{theorem}\label{sevcompt} Assume that the Gorenstein sequence $H$
satisfies $H=(1,4,7,h,b,\ldots ,4,1)$, of socle degree $j\ge 6$, where
$h,b$ are arbitrary integers satisfying the necessary restrictions of
Lemma~\ref{nonempty2a}.
\renewcommand{\theenumi}{\roman{enumi}}
\begin{enumerate}[i.]
\item\label{sevcompti} the dimension of the tangent space
$\mathcal T_I$ on $\PGOR(H)$ to a general element $I$ of\linebreak
$\mathfrak{ C}(H)\subset
\PGOR(H)$ satisfies,
\begin{equation}\label{tangspe2}
 \dim_K \mathcal{T}_I=\dim \mathfrak{C}(H) +1+\nu_{j-1}(J)
\end{equation}
where $J$ is a generic element of $\PGOR(H'), H'=(1,3,6,h-1,b-1,
\ldots, h-1,6,3,1)$.
\item\label{sevcomptii} When $j\ge 6$, the Zariski closure
$\overline{\mathfrak{C}(H)}$ is a
generically smooth irreducible component of $\PGOR(H)$ when, equivalently
\begin{enumerate}
\item\label{sevcomptiia} $\nu_{j-1}(J)=0$ for $J$ generic in $\PGOR(H')$;
\item\label{sevcomptiib} a generic $J\in\PGOR(H')$ has no degree-4
relations;
\item\label{sevcomptiic} $3h-b-17\ge 0$.
\end{enumerate}
\end{enumerate}
\end{theorem}
\begin{proof}  Here \eqref{sevcompti} follows immediately from Theorem
\ref{sevcomp}
\eqref{sevcompi},\eqref{sevcompii}. This shows
\eqref{sevcomptiia}; by the symmetry of the minimal resolution of $J$,
\eqref{sevcomptiia} is equivalent to \eqref{sevcomptiib}. The third difference
satisfies $(\Delta^3 H')_4=17+b-3h$, and under the assumption $j\ge 6$, it
gives, when
positive, the number of degree-4 relations --- the linear relations among those
generators of $J$ having degree 3; when 0 or negative there are no such
relations. This
completes the proof of the equivalence of \eqref{sevcomptiib} and
\eqref{sevcomptiic}.
\end{proof}\par
We now show that there are monomial ideals in $R'=K[x,y,z]$, having certain
Hilbert functions $T'$ and
having a small number of generators.  This prepares a key step for
Theorem
\ref{sevcomp2}. We consider Hilbert functions of the form
$T'=(1,3,3,\ldots , 2_a,\ldots ,1_c,\ldots ,0,\ldots )$ where degree $a$ is
the first degree in which
$T'_a<3$, and $c$ is the first degree $c\ge 3$ in which $T'_i\le 1$, and
$d$ is the first positive degree in
which $T'_i=0$: we allow equalities among $a,c,d$, so if $a=c=4,d=5$,
$T'=(1,3,3,3,1,0,\ldots )$. The
following result is easy to verify.
\begin{lemma}\label{sevcomp2AB} \begin{enumerate}[i.]
\item\label{sevcomp2A} The Artinian algebra $A=R'/J_{a,c,d},
J_{a,c,d}=(xy,xz,yz,x^a,y^c,z^d), 3\le a\le c\le
d$ has Hilbert function
$T'(a,c,d)=(1,3,3,\ldots ,2_a,\ldots ,1_c,\ldots ,0_d,\ldots )$  in the
sense above. \par
\item\label{sevcomp2B} The Artinian algebra $A=R'/K_{a,c},
K_{a,c}=(x^2,xy,z^2,x^{a-1}z,y^c),3\le a\le c$ has
Hilbert function $T'(a,c)=(1,3,3,2,\ldots ,1_a,\ldots ,0_c,\ldots )$.
\end{enumerate}
\end{lemma}
\begin{corollary}{\sc{Artinian Gorenstein algebras with related minimal
resolution.}}\label{relGor}
\begin{enumerate}[i.]
\item (P. Maroscia)
\cite{Mar,GMR,IK,MiN} Let
$s=\sum_{i\ge 0} T'(a,c,d)_i$, or
$\sum _{i\ge0}T'_{a,c}$, respectively. Then there are smooth degree-$s$
punctual schemes $\Z=\Z (a,c,d)\subset
\mathbb P^3$ or $\Z=\Z (a,c)\subset \mathbb P^3$, respectively, whose
coordinate rings have the same minimal
resolutions as the Artinian algebras defined by $J_{a,c,d}$ or $K_{a,c}$,
respectively.
\par
\item  (M. Boij)
\cite{Bj1})\label{goodcomp}
Furthermore, let $j\ge 2c$, or $j\ge 2b$, respectively, and let
$A=A(a,c,d,j,F)$ or $A=A(a,c,j,F)$,
respectively, denote a general enough GA quotient of $\mathcal O_\Z,
\Z=\Z(a,c,d)$ or $\Z=\Z(a,c)$
having socle degree $j$, defined by
$A=R/\Ann (F), F\in (I_\Z)_j^\perp $. The minimal resolution of
$A$ agrees with that of the corresponding coordinate ring $\mathcal O_\Z$
in degrees up to $j/2$.
\end{enumerate}
\end{corollary}
\begin{proof} P. Maroscia's well known result deforms a given monomial
ideal defining an Artinian
algebra to a graded ideal defining a smooth punctual scheme $\Z$, and
having the same minimal
resolution. M.~Boij showed that a general enough GA quotient of $\Z$ has a
related minimal
resolution.
\end{proof}
\begin{theorem}\label{sevcomp2}{\sc Families $\PGOR(H)$ with several
components.}
\begin{enumerate}[i.]
\item\label{sevcomp2A} Assume that $H$ is a Gorenstein sequence of socle
degree $j\ge 6$ satisfying
\eqref{genhfe}, namely $H=(1,4,7,h,b,\ldots )$ and that $h\le 10$. Then
there is a GA quotient of
the coordinate ring of a smooth punctual scheme $\Z$ having Hilbert
function $H$, and
$H=\Sym(H_\Z,j)$.
\item\label{sevcomp2B} Assume further that
$3h-b-17\ge 0$ and $8\le h\le 10$. Then $\PGOR(H)$ has at least two
irreducible components, the
component
$\overline{\mathfrak {C}}$, and a component containing
suitable GA algebras
$A=A(a,c,d,j,F)$ or $A=A(a,c,j,F)$,
respectively that are quotients of the coordinate ring of smooth punctual
schemes.
\end{enumerate}
\end{theorem}
\begin{proof} Assume that $H=(1,4,7,h,b,\ldots ,1)$ has socle
degree $j\ge 6$ and let $T'=\Delta H_{\le
j/2}$. By Theorem
\ref{nonempty2},
$T'$ is an $O$-sequence; since $h \le 10$ Lemma \ref{nonempty2a} implies
$T'$ satisfies
$T'=(1,3,3,h-7,b-h,\ldots )$, with
$h-7\le 3$, with $T'_{\le j/2}$ nonnegative, and nonincreasing after
degree $0$ to $1$. Thus
$T'=
T'(a,c,d)$ or
$T'=T'(a',c')$ for suitable $(a,c,d)$ or $(a',c')$. Lemma \ref{sevcomp2AB}
and Corollary
\ref{relGor}\eqref{goodcomp}
imply that there is a Artinian Gorenstein algebra $A=R/I$ of Hilbert
function $H$, such that the
beginning of its minimal resolution is that of $R'/J(a,b,c)$ or
$R'/K(a,b)$. In particular $I_2$
has at most two linear relationss. Since one cannot specialize from a GA
algebra
$A=R/I\in
\overline{\mathfrak {C}(H)}$  where
$I_2$ has three linear relations, to a GA algebra $A=A(a,c,d,j,F)$ or
$A(a,c,j,F)$ where $I_2$ has at most
two linear relations, the claim of the theorem follows.
\end{proof}

\smallskip \par \noindent
{\bf Acknowledgements}\par \smallskip The authors thank Carol Chang, Dale
Cutkosky, Juan Migliore,  Hal
Schenck, and Jerzy Weyman for helpful discussions; we are also appreciative
to Joe Harris who loaned us a
copy of Y.~A.~Lee's unpublished senior thesis, completed under his
direction. The authors thank
the referees for helpful comments, in particular one leading to our use of
Corollary \ref{oseq}.
\medskip
\bibliographystyle{amsalpha}

\end{document}